\documentclass[12pt,letterpaper]{article}
\usepackage[utf8]{inputenc}
\usepackage{amsmath}
\usepackage{amsfonts}
\usepackage{amssymb}
\usepackage{amsthm}  
\usepackage[top=1in, bottom=1in, left=1in, right=1in]{geometry}
\usepackage{mathdots}
\usepackage{bbold}
\usepackage{titlesec}
\usepackage{hyperref}
\usepackage[bottom]{footmisc}
\usepackage{mathrsfs}
\usepackage{graphicx}
\usepackage{verbatim}
\usepackage{algpseudocode}
\usepackage{textcomp}
\usepackage{xcolor}
\usepackage{qtree} 
\titleformat{\section}{\normalsize\scshape\center}{\thesection}{1em}{}
\titleformat{\subsection}{\normalsize\scshape\center}{\thesubsection}{1em}{}
\makeatletter
\renewcommand\@makefntext[1]{%
    \parindent 1em%
    \@thefnmark.~#1}
\makeatother

\makeatletter
\def\moverlay{\mathpalette\mov@rlay}
\def\mov@rlay#1#2{\leavevmode\vtop{%
   \baselineskip\z@skip \lineskiplimit-\maxdimen
   \ialign{\hfil$\m@th#1##$\hfil\cr#2\crcr}}}
\newcommand{\charfusion}[3][\mathord]{
    #1{\ifx#1\mathop\vphantom{#2}\fi
        \mathpalette\mov@rlay{#2\cr#3}
      }
    \ifx#1\mathop\expandafter\displaylimits\fi}
\makeatother

\newcommand{\id}{\textnormal{id}}

\newtheorem{definition}{Definition}[section]
\newtheorem{theorem}[definition]{Theorem}
\newtheorem{proposition}[definition]{Proposition}
\newtheorem{lemma}[definition]{Lemma}
\newtheorem{corollary}[definition]{Corollary}
\newtheorem{conjecture}[definition]{Conjecture}

\newtheorem{question}[definition]{Question}

\let\OLDthebibliography\thebibliography
\renewcommand\thebibliography[1]{
  \OLDthebibliography{#1}
  \setlength{\parskip}{0pt}
  \setlength{\itemsep}{0pt plus 0.3ex}
}

\begin{document}
\title{\uppercase{\textbf{\normalsize The combinatorics of a tree-like functional equation for connected chord diagrams}}}
\author{\small{\textsc{Lukas Nabergall}}}
\date{\small{\textsc{\today}}}
\maketitle

\allowdisplaybreaks

\begin{abstract}
We build on recent work of Yeats, Courtiel, and others involving connected chord diagrams. We first derive from a Hopf-algebraic foundation a class of tree-like functional equations and prove that they are solved by weighted generating functions of two different subsets of weighted connected chord diagrams: arbitrary diagrams and diagrams forbidding so-called top cycle subdiagrams. These equations generalize the classic specification for increasing ordered trees and their solution uses a novel decomposition, simplifying and generalizing previous results. The resulting tree perspective on chord diagrams leads to new enumerative insights through the study of novel diagram classes. We present a recursive bijection between connected top-cycle-free diagrams with $n$ chords and triangulations of a disk with $n+1$ vertices, thereby counting the former. This connects to combinatorial maps, Catalan intervals, and uniquely sorted permutations, leading to new conjectured bijective relationships between diagram classes defined by forbidding graphical subdiagrams and imposing connectedness properties and a variety of other rich combinatorial objects. We conclude by exhibiting and studying a direct bijection between diagrams of size $n$ with a single terminal chord and diagrams of size $n-1$. 

\end{abstract}

\section{Introduction}

A chord diagram of size $n$ is a perfect matching of $\{1, 2, \ldots, 2n\}$. Such objects have also been called matchings \cite{Jelinek2005, Chen2007}, linked diagrams \cite{Stein1978a}, complete pairings \cite{Stein1978a}, and interval systems \cite{Davies2020}, and can be viewed as set partitions with every block of size 2; see Figure~\ref{chord_diagram_fig}. Here chord diagrams are rooted at the chord containing 1, but the term has also been used to refer to the unrooted object obtained after modding out by cyclic permutations of $[2n]$. For brevity, we will use the shortened `diagram' to refer to chord diagrams. There are $(2n - 1)!!$ diagrams of size $n$. Other objects counted by double factorials include double occurrence words, fixed-point-free involutions, increasing trees, and Stirling permutations. In the literature, chord diagrams seem to have been first studied by Touchard \cite{Touchard1952}. Since then there has been much focus on the enumeration of various subclasses of diagrams and their statistics (e.g. \cite{Riordan1975, Stein1978a, Stein1978b, Nijenhuis1979, Flajolet2000, Pilaud2014, Courtiel2017b}). One of the most prominent and natural types of diagrams studied are connected diagrams, which are those diagrams $C$ for which there is no proper interval of $[2n]$ that is the ground set of a subdiagram of $C$. Connectivity can also be defined via the intersection graph $G(C)$ of a diagram $C$, the directed graph on the chords of $C$ formed by adding edge $(c, c')$ if $c'$ crosses $c$ on the right; $C$ is connected if and only if its intersection graph is weakly connected. Undirected intersection graphs of diagrams, known as circle graphs, have been well studied as pure graph objects (e.g. \cite{Naji1985, Davies2020}) and in particular play a major role in vertex minor theory \cite{Bouchet1994}. 

Outside of enumerative combinatorics and graph theory, chord diagrams have also appeared in diverse areas such as knot theory \cite{Bollobas2000}, bioinformatics \cite{Hofacker1998} and, most relevantly for the present purposes, physics. In 2013, Nicolas Marie and Karen Yeats \cite{Marie2013} solved a certain Dyson-Schwinger equation from quantum field theory. This equation has a recursive form similar to standard functional equations for rooted trees. The solution of Marie and Yeats came as a series expansion indexed by connected diagrams and weighted by coefficients indexed by certain novel parameters of those diagrams. These parameters are the indices and count of the terminal chords of a diagram $C$ under a total order on the chords of $C$, the intersection order; a chord $c \in C$ is terminal if it has no outgoing edges in $G(C)$, that is, if there are no chords crossing it to the right. The intersection order is a natural total order, distinct from the standard ordering of the chords by their first endpoints or sources, that extends the partial order on the chords defined by reachability in the intersection graph. See Section~\ref{chord_diagram_sec} for more details. Following Marie and Yeats, Courtiel and Yeats \cite{Courtiel2017a, Courtiel2019} further studied the expansion over connected diagrams and the associated terminal chords, in particular obtaining probabilistic and asymptotic information about the distribution of the latter. In a different direction, Courtiel, Yeats, and Zeilberger \cite{Courtiel2017b} linked these objects and parameters to combinatorial maps, another well-studied structure. We will further describe related work throughout the paper. 

We aim to bring these new equations, parameters, and associated results further into the enumerative combinatorial context and, through this, obtain new insights. In particular, much of this paper can be viewed as illuminating the tree-like structure of chord diagrams and its interplay with notions of connectedness, as well as relating these insights to other combinatorial objects. We begin in Section~\ref{hopf_sec} by showing how generalizations of the Dyson-Schwinger equation considered by Marie and Yeats arise from Hopf subalgebras of the Connes-Kreimer Hopf algebra of rooted trees via work by Lo\"{i}c Foissy \cite{Foissy2007, Foissy2010, Foissy2014}. These equations have the form
\begin{align}
\label{tree_like_eq}
G(x, y) = xL(\phi(G(x, y))),
\end{align}
where $\phi$ is a formal power series with nonzero constant term. The map $L$ is a type of linear map on polynomials in $y$ from Hochschild cohomology theory known as a {\em Hochschild 1-cocycle}, which is defined by the identity
\begin{align*}
\Delta \circ L = (\textnormal{id} \otimes L)\circ \Delta + L \otimes 1,
\end{align*}
where $\Delta$ is the coproduct of a coalgebra, $\textnormal{id}$ is the identity map, and $1$ is the identity of the underlying ring. There are two such 1-cocyles corresponding to the two coalgebras on the ring of one-variable polynomials, the binomial coalgebra and the divided power coalgebra, leading to two forms of equation (\ref{tree_like_eq}). Since in the binomial case $L$ generalizes the integral operator, (\ref{tree_like_eq}) generalizes the classic functional equation for the exponential generating function of increasing trees (see e.g. \cite{Bergeron1992}), leading to the term {\em tree-like} to describe such equations (matching the traditional generalized meaning of the term \cite{Bergeron1998}). 

Marie and Yeats \cite{Marie2013} solved the tree-like equation with $L$ a binomial 1-cocycle and $\phi(z) = 1/(1 - z)$. Their solution method proceeded in two steps. First they applied a decomposition on connected diagrams to prove that the solution satisfies a certain recurrence, in the form of a differential equation, which allows them to reduce the problem to verifying the linear term in $y$. Then they inductively expanded the reduced form of (\ref{tree_like_eq}) to obtain a second recurrence characterizing the solution and proved that it holds by passing through a bijection between connected diagrams and a recursively-defined class of labeled binary trees. Their proof method is quite technical and therefore difficult to understand intuitively and generalize. In Sections~\ref{func_to_identity_sec} through \ref{func_eq_sol_sec} we present a simpler, more direct proof of our main result, Theorem~\ref{equation_sol_thm}, generalizing the work of Marie and Yeats, that the binomial and divided power tree-like equations are solved by weighted generating functions for sets of connected chord diagrams. Along the way, in Section~\ref{chord_diagram_sec}, we define chord diagrams and their relevant properties and parameters, including the intersection order, terminal chords, and $k$-terminality, a property that we later show can be viewed as a stronger form of $k$-connectivity. Our proof of Theorem~\ref{equation_sol_thm} works entirely at the level of chord diagrams and is based on a relatively simple decomposition of a connected diagram into ``shuffled" connected subdiagrams. This decomposition naturally involves distinguished 1-terminal subdiagrams, beginning to highlight the importance of 1-terminality; in particular, these 1-terminal subdiagrams essentially reveal the underlying tree-like structure of chord diagrams. 

In Section~\ref{diff_eq_sec} we turn towards further explaining the appearance of chord diagrams and how they relate to the choice of linear map $L$ in the tree-like equation by proving that, with $\phi(z) = 1/(1 - z)$, $G$ satisfies the differential equation from \cite{Marie2013} if and only if $L$ is a binomial 1-cocycle. This differential equation is revealing because it is exactly a two variable generalization of the well-known functional equation
\begin{align*}
C(x) - x = C(x)\left(2x{d \over dx} - 1\right)C(x)
\end{align*} 
for the generating function $C(x)$ of connected diagrams. The ``if" direction of this result is more straightforward than the other direction and was already proved as a key ingredient in \cite{Marie2013}.  

The solution to the divided power tree-like equation is indexed by a subset of connected diagrams, namely, those that forbid so-called top cycle subdiagrams, one of only two types of chord diagrams whose undirected intersection graph is isomorphic to a cycle. In Section~\ref{comb_map_triangle_sec} we describe the bijective relationship between connected top-cycle-free diagrams and planar bridgeless maps discovered by Courtiel, Yeats, and Zeilberger \cite{Courtiel2017b}. This left open the goal, further motivated by our generating function results, of obtaining explicit formulas counting such diagrams. We derive such a formula counting the number of connected top-cycle-free diagrams of size $n$ by describing an explicit, recursive bijection to triangulations of a disk and applying the work of Brown \cite{Brown1964} enumerating these classic objects. This formula is in fact further refined by the index of the first terminal chord in the intersection order, which turns out to correspond to one less than the number of exterior (or boundary) vertices of the triangulation. Alongside prior work by Jel\'{i}nek \cite{Jelinek2005}, these results motivate considering other classes of diagrams determined by forbidding a fixed set of graphically-defined subdiagrams, in the vein of the large body of prior work on analogous classes of graphs as well as pattern-avoiding permutations. In Sections \ref{interval_sorted_sec} and \ref{forbid_subdiagram_sec} we discuss relationships, both conjectured and proven, between top-cycle-free, triangle-free, tree, chordal, bipartite, and bottom-cycle-free diagrams, where a bottom cycle is the second diagram realizer of a cycle, and intervals of Catalan lattices, uniquely sorted permutations, and a variety of combinatorially-significant number sequences found on the OEIS \cite{OEIS2021}. These relationships in particular appear when the diagrams are required to be connected or 1-terminal, further underscoring the role of connectivity notions. Resolving the number of conjectures introduced here is a significant line of inquiry for future work. 

Courtiel and Yeats \cite{Courtiel2017a} proved that there are $(2n-3)!!$ 1-terminal diagrams of size $n$, indicating that such diagrams should be in one-to-one correspondence with diagrams of size $n-1$. In the final section of this paper we present a new simple formulation of an unpublished bijection between these two sets first discovered by Yeats. This map $\psi$ is easy to understand structurally. It also has a straightforward relationship with $k$-terminality; $\psi$ induces a bijection between $k$-terminal diagrams of size $n$ and $(k-1)$-terminal diagrams of size $n-1$. Together with a characterization of 1-terminal top-cycle-free diagrams, it follows from this that such diagrams are in bijection with noncrossing diagrams of size $n-1$, so they are counted by the Catalan numbers. We also show that $\psi$ interfaces equally well with nonnesting diagrams, another classic Catalan object. After briefly discussing the action of $\psi$ on other double factorial objects, we conclude by proving that applying $\psi$ necessarily decreases the connectivity of a diagram by at least 1. This result is a step in the direction towards a complete understanding of the relationship between 1-terminality, the map $\psi$, and connectivity.







\section{Tree-like equations from the Connes-Kreimer Hopf algebra}
\label{hopf_sec}

First appearing under the guise of the Butcher group in numerical analysis and independently introduced by Kreimer \cite{Kreimer1998} in the context of renormalization in perturbative quantum field theory, the Connes-Kreimer Hopf algebra of rooted trees $\mathcal{H}_{CK}$ is the free associative commutative algebra freely generated over a field $K$ of characteristic zero by the set of rooted trees. As implied, the product is given on the basis of forests of rooted trees by concatenation while the coproduct is defined by setting
\begin{align*}
\Delta(t) = \sum_{\substack{C \subseteq V(t) \\ C \textrm{ antichain}}}\left(\prod_{v \in C}t_{v}\right)  \otimes \left(t \setminus \prod_{v \in C}t_{v}\right)
\end{align*}
for a rooted tree $t$, where $t_{v}$ is the subtree of $t$ rooted at $v$, and extending it as an algebra homomorphism to all of $\mathcal{H}_{CK}$; note that here we take $t \setminus t = 1$. From a pure algebraic perspective, the Connes-Kreimer Hopf algebra is important because it possesses a certain {\em universal property} unique up to isomorphism among Hopf algebras:

\begin{theorem}[Connes-Kreimer {\cite[Theorem 2]{Connes1998}}]
\label{univ_prop_thm}
Let $A$ be an associative commutative algebra over $K$ and $L: A \to A$ be a linear map. Then there exists a unique algebra homomorphism $\rho_{L}: \mathcal{H}_{CK} \to A$ such that $\rho_{L} \circ B_{+} = L \circ \rho_{L}$. Furthermore, if $A$ is a bialgebra and $L$ is a Hochschild 1-cocycle then $\rho_{L}$ is a bialgebra homomorphism, and if $A$ is also a Hopf algebra then $\rho_{L}$ is a Hopf algebra homomorphism. 
\end{theorem}

There has been considerable interest in understanding Hopf subalgebras of $\mathcal{H}_{CK}$. In a series of papers \cite{Foissy2007, Foissy2010, Foissy2014}, Lo\"{i}c Foissy examined subalgebras of $\mathcal{H}_{CK}$ generated by a family of recursive equations, so-called combinatorial Dyson-Schwinger equations, of the form
\begin{align}
\label{tree_eq}
T(x) = xB_{+}(\phi(T(x)))
\end{align}
for $\phi(z) \in K[[z]]$ with $\phi(0) = 1$. The unique solution to this equation is a formal power series $T(x)$ whose coefficients lie in $\mathcal{H}_{CK}$. Writing $t_{n} = [x^{n}]T(x)$, Foissy characterized when the subalgebra $A = K[t_{1}, t_{2}, \ldots]$ of $\mathcal{H}_{CK}$ is Hopf. 

\begin{theorem}[Foissy \cite{Foissy2007}]
\label{hopf_subalgebra_thm}
$A$ is a Hopf subalgebra if and only if $\phi(z) = (1 + abz)^{-1/b}$ for some $a, b \in K$ with $b \neq 0$ or $\phi(z) = e^{az}$. 
\end{theorem}

We are interested in equations which arise from (\ref{tree_eq}) by applying the universal property to the polynomial algebra $K[y]$ and a linear map $L: K[y] \rightarrow K[y]$. Applying the algebra homomorphism $\rho_{L}$ guaranteed by Theorem~\ref{univ_prop_thm} to both sides of (\ref{tree_eq}), we get the bivariate {\em tree-like equation}
\begin{align}
\label{gen_func_eq1}
G(x, y) = xL(\phi(G(x, y))),
\end{align} 
where $G(x, y) = \rho_{L}(T(x))$. The maps $L$ and $\rho_{L}$ act on the coefficients in $x$ of $\phi(G(x, y))$ and $T(x)$ term by term additively; since $L$ sends polynomials to polynomials this equation has an inductively specified solution in $K[y][[x]]$, so it is well-formed. In the physics setting, $\rho_{L}$ corresponds to the Feynman rules which map each Feynman graph to its associated Feynman integral (for details see e.g. \cite{Panzer2011}). We will be most interested in equation (\ref{gen_func_eq1}) with $\phi$ set to generate a Hopf subalgebra of $\mathcal{H}_{CK}$ via Theorem~\ref{hopf_subalgebra_thm}, but will work in the more general setting with $\phi$ an arbitrary formal power series with constant term 1. 

In order to get meaningful combinatorial solutions to (\ref{gen_func_eq1}), it is clearly necessary to restrict $L$ to some specific class of linear maps. With that in mind, the universal property points the way towards which classes of maps would be of most interest: Hochschild 1-cocycle operators arising from coalgebra structures on $K[y]$. There are two graded coalgebras on one-variable polynomials classically studied in the literature, namely, the binomial coalgebra and the divided power coalgebra. For the former, the coproduct is defined by setting
\begin{align*}
\Delta(y^{n}) = \sum_{k=0}^{n}{n \choose k}y^{k} \otimes y^{n-k}.
\end{align*}
Combining this with the polynomial algebra on $K[y]$ with the usual product, we in fact get a Hopf algebra. The following lemma describes Hochschild 1-cocycles in the binomial coalgebra. 

\begin{lemma}
\label{bin_1cocycle_lem}
If $L$ is a 1-cocycle operator for the binomial coalgebra on $K[y]$, then 
\begin{align*}
L(y^{n}) = \int_{0}^{y}F{\left({d \over dt}\right)}t^{n}dt
\end{align*}
for some power series $F(z) = \sum_{i \geqslant 0}f_{i}z^{i}$ in $K[[z]]$. 
\end{lemma}
\begin{proof}
Writing $c_{m,n} = [y^{n-m}]L(y^{n})$ for $m \leqslant n$, define
\begin{align*}
L_{m}(y^{n}) = \begin{cases}
c_{m, n}y^{n-m} &\textrm{if } n \geqslant m \\
0 &\textrm{else.}
\end{cases}
\end{align*}
Fix $m \leqslant n$ and note that $L_{m}$ is a 1-cocycle by linearity. Then
\begin{align*}
((\textnormal{id} \otimes L_{m})\circ \Delta + L_{m} \otimes 1)(y^{n}) &= (\textnormal{id} \otimes L_{m})\left(\sum_{k=0}^{n}{n \choose k}y^{k} \otimes y^{n-k}\right) + L_{m}(y^{n}) \otimes 1 \\
&= \sum_{k=0}^{n}c_{m,n-k}{n \choose k}y^{k} \otimes y^{n-m-k} + c_{m,n}y^{n-m} \otimes 1,
\end{align*}
while 
\begin{align*}
(\Delta \circ L_{m})(y^{n}) = c_{m,n}\Delta(y^{n-m}) = \sum_{i=0}^{n-m}c_{m,n}{n-m \choose i}y^{i} \otimes y^{n-m-i}.
\end{align*}
Applying the 1-cocycle property and comparing terms, we see that $m \geqslant -1$ since otherwise $x^{n+1} \otimes x^{-m-1}$ appears in the ladder but not the former. Furthermore, 
\begin{align*}
c_{m,n}{n - m \choose k} = \begin{cases}
c_{m,n-k}{n \choose k} &\textrm{if } k \neq n - m \\ 
c_{m,n-k}{n \choose k} + c_{m,n} &\textrm{if } k = n - m.
\end{cases}
\end{align*}
It follows that $c_{m,m} = 0$ and, for $m < n$, $c_{m,n} = {n! \over (n-m)!(m+1)!}c_{m,m+1}$. One can then readily check that we get the desired expression for $L$ by setting $f_{m+1} = c_{m,m+1}/(m+1)!$. 
\end{proof}

Although we proved it for completeness, this result is well known; e.g. Panzer \cite{Panzer2011} obtained an equivalent algebraic characterization. For the divided power coalgebra, the coproduct is defined by setting 
\begin{align*}
\Delta(y^{n}) = \sum_{k=0}^{n}y^{k} \otimes y^{n-k}.
\end{align*}
This also gives a Hopf algebra on $K[y]$, but the compatible algebra structure instead has the product $y^{i} \cdot y^{j} = {i + j \choose i}y^{i+j}$; nevertheless it is in fact isomorphic as a Hopf algebra to the binomial Hopf algebra via a scaling of coefficients. With that said, we will always work with the standard product on $K[y]$, meaning that in this case only the first statement of the universal property will apply; we will later see that interesting combinatorics arise regardless. A similar formula holds for Hochschild 1-cocycles in the divided power coalgebra with the integral replaced with a degree raising operator and the derivative with a degree lowering operator ${\delta \over \delta y}$---its proof can be easily constructed by adapting the proof of Lemma~\ref{bin_1cocycle_lem} so we omit it. 

\begin{lemma}
\label{power_1cocycle_lem}
If $L$ is a 1-cocycle operator for the divided power coalgebra on $K[y]$, then 
\begin{align*}
L(y^{n}) = yF{\left({\delta \over \delta y}\right)}y^{n}
\end{align*}
for some power series $F(z) = \sum_{i \geqslant 0}f_{i}z^{i}$ in $K[[z]]$, where ${\delta \over \delta y}y^{n} = y^{n-1}$ if $n > 0$ and $0$ otherwise. 
\end{lemma} 

We will write $L_{bin}$ and $L_{div}$ for 1-cocycles of the binomial and divided power coalgebras, respectively, with the underlying power series $F$ implicitly carried along. Note that for both of these operators $\deg L(y^{n}) \leqslant n+1$ and, in particular, this bound is obtained if and only if $f_{0} \neq 0$.

\section{From functional equation to combinatorial identity}
\label{func_to_identity_sec}

Over the next three sections, we will prove that 
\begin{align}
\label{bin_func_eq}
G(x, y) = xL_{bin}(\phi(G(x, y)))
\end{align}
and 
\begin{align}
\label{power_func_eq}
G(x, y) = xL_{div}(\phi(G(x, y)))
\end{align}
are solved by certain weighted generating functions over connected weighted chord diagrams. Our proof strategy is as follows: 1) expand the functional equations into recurrences which characterize the solution, 2) apply induction to turn the recurrence into an equivalent generic combinatorial identity, then for each case, 3) set up a map underlying the corresponding identity based on a decomposition of a weighted chord diagram into weighted subdiagrams, 4) exhibit its inverse, proving that it is a bijection and therefore the identity holds, as required. 

We begin with the expansion into a recurrence, starting with the generic tree-like equation
\begin{align}
\label{gen_func_eq2}
G(x, y) = xL(\phi(G(x, y)))
\end{align}
for $\phi(z) = \sum_{k \geqslant 0}\phi_{k}z^{k}$ with $\phi_{0} = 1$ and a linear map $L$ on the binomial algebra $K[y]$. Since all the cases we will be concerned with in this paper involve polynomial operators that raise the degree in $y$ by at most 1, we assume that $L$ is such an operator; in particular, $\deg L(y^{n}) \leqslant n + 1$ for all $n \geqslant 0$. 
Since $L$ maps polynomials to polynomials of degree at most one higher and $G$ clearly has no constant term, $G(x, y) = \sum_{i \geqslant 1}h_{i}(y)x^{i}$ for some polynomials $h_{i}(y)$ of degree at most $i$. Expanding the composition $\phi(G(x, y))$ of (\ref{gen_func_eq2}) and applying linearity of $L$ gives
\begin{align*}
G(x, y) &= xL\left(\sum_{k \geqslant 0}\phi_{k}\left(\sum_{i \geqslant 1}h_{i}(y)x^{i}\right)^{k}\right) \\
&= x\phi_{0}L(1) + xL\Bigg(\sum_{k \geqslant 1}\sum_{n \geqslant 1}\sum_{\substack{n_{1} + \cdots + n_{k} = n \\ n_{\ell} \geqslant 1}}\phi_{k}h_{n_{1}}(y)\cdots h_{n_{k}}(y)x^{n}\Bigg) \\
&= xL(1) + x\sum_{k \geqslant 1}\sum_{n \geqslant 1}\sum_{\substack{n_{1} + \cdots + n_{k} = n \\ n_{\ell} \geqslant 1}}\phi_{k}L(h_{n_{1}}(y)\cdots h_{n_{k}}(y))x^{n}
\end{align*}
Comparing coefficients, we obtain the recurrence
\begin{align}
\label{recurrence_eq}
h_{n+1}(y) = L(1)\mathbb{1}_{n=0} + \sum_{k=1}^{n}\sum_{\substack{n_{1} + \cdots + n_{k} = n \\ n_{\ell} \geqslant 1}}\phi_{k}L(h_{n_{1}}(y)\cdots h_{n_{k}}(y)),
\end{align}
where $\mathbb{1}_{n=0} = 1$ if $n = 0$ and $0$ otherwise. Since clearly all of the above steps are reversible, this recurrence uniquely specifies the solution to (\ref{gen_func_eq2}). 

Now we translate this recurrence into an equivalent combinatorial identity. To do this without specifying beforehand the linear map $L$ we will write
\begin{align}
\label{gen_map_eq}
L(y^{n}) = a_{n}\sum_{i=0}^{n+1}f_{n,i}b_{i}y^{i}
\end{align}
for the action of $L$ on the standard basis of $K[y]$ and suppose that the solution has the form
\begin{align}
\label{gen_solution_eq}
h_{n}(y) = \sum_{\substack{A \in \mathcal{A} \\ |A| = n}}f_{A}\phi_{A}{L(y^{m(A)-1}) \over a_{m(A)-1}},
\end{align}
that is, $G(x, y)$ is a kind of generating function of a set of $\phi$-weighted combinatorial objects $\mathcal{A}$ counted by size $n(A) = |A|$ in the $x$ variable and one less than another integer parameter $m(A) \leqslant n(A)$ in the $y$ variable. This generating function is weighted by terms $f_{A}$ determined by the coefficients $f_{n,i}$ defining $L$ and the multiplicative inverse of a constant $a_{m(A)-1} \in K$ determined only by $m(A)$. The $\phi$-weighting appears as an additional term $\phi_{A}$ which is a product of the weights of $A$. We will comment further on this kind of generating function when specializing the argument to solve (\ref{bin_func_eq}) and (\ref{power_func_eq}). 

To perform the translation, we apply induction on $n$. Substituting (\ref{gen_map_eq}) into (\ref{gen_solution_eq}), we get
\begin{align*}
h_{n}(y) = \sum_{\substack{A \in \mathcal{A} \\ |A| = n}}f_{A}\phi_{A}\sum_{i=0}^{m(A)}f_{m(A)-1,i}b_{i}y^{i} = \sum_{i=0}^{n}\Bigg(\sum_{\substack{A \in \mathcal{A} \\ |A| = n \\ m(A) \geqslant i}}f_{m(A)-1, i}f_{A}\phi_{A}\Bigg)b_{i}y^{i},
\end{align*}
thereby expressing $h_{n}(y)$ in standard polynomial form. Expanding each term in the nested series on the right hand side of (\ref{recurrence_eq}), we inductively have
\begin{align*}
& \phi_{k}L(h_{n_{1}}(y)\cdots h_{n_{k}}(y)) \\ 
&\quad = \phi_{k}L\Bigg(\sum_{i_{1}=0}^{n_{1}}\Bigg(\sum_{\substack{A_{1} \in \mathcal{A} \\ |A_{1}| = n_{1} \\ m(A_{1}) \geqslant i_{1}}}f_{m(A_{1})-1, i_{1}}f_{A_{1}}\phi_{A_{1}}\Bigg)b_{i_{1}}y^{i_{1}} \cdots \sum_{i_{k}=0}^{n_{k}}\Bigg(\sum_{\substack{A_{k} \in \mathcal{A} \\ |A_{k}| = n_{k} \\ m(A_{k}) \geqslant i_{k}}}f_{m(A_{k})-1, i_{k}}f_{A_{k}}\phi_{A_{k}}\Bigg)b_{i_{k}}y^{i_{k}}\Bigg) \\
&\quad = \phi_{k}\sum_{m=0}^{n}L(y^{m})\sum_{\substack{i_{1} + \cdots + i_{k} = m \\ 0 \leqslant i_{\ell} \leqslant n_{\ell}}}b_{i_{1}}\cdots b_{i_{k}}\Bigg(\sum_{\substack{A_{1} \in \mathcal{A} \\ |A_{1}| = n_{1} \\ m(A_{1}) \geqslant i_{1}}}f_{m(A_{1})-1, i_{1}}f_{A_{1}}\phi_{A_{1}}\Bigg) \dots \Bigg(\sum_{\substack{A_{k} \in \mathcal{A} \\ |A_{k}| = n_{k} \\ m(A_{k}) \geqslant i_{k}}}f_{m(A_{k})-1, i_{k}}f_{A_{k}}\phi_{A_{k}}\Bigg) \\
&\quad = \phi_{k}\sum_{m=0}^{n}\sum_{i=0}^{m+1}f_{m,i}b_{i}y^{i}\sum_{\substack{i_{1} + \cdots + i_{k} = m \\ 0 \leqslant i_{\ell} \leqslant n_{\ell}}}a_{m}b_{i_{1}}\cdots b_{i_{k}} \\ 
&\qquad\qquad\qquad\qquad \times \Bigg(\sum_{\substack{A_{1} \in \mathcal{A} \\ |A_{1}| = n_{1} \\ m(A_{1}) \geqslant i_{1}}}f_{m(A_{1})-1, i_{1}}f_{A_{1}}\phi_{A_{1}}\Bigg) \dots \Bigg(\sum_{\substack{A_{k} \in \mathcal{A} \\ |A_{k}| = n_{k} \\ m(A_{k}) \geqslant i_{k}}}f_{m(A_{k})-1, i_{k}}f_{A_{k}}\phi_{A_{k}}\Bigg) \\
&\quad = \sum_{i=0}^{n+1}b_{i}y^{i}\sum_{m=\max\{i-1,0\}}^{n}f_{m,i}\phi_{k}\sum_{\substack{i_{1} + \cdots + i_{k} = m \\ 0 \leqslant i_{\ell} \leqslant n_{\ell}}}a_{m}b_{i_{1}}\cdots b_{i_{k}} \\ 
&\qquad\qquad\qquad\qquad \times \Bigg(\sum_{\substack{A_{1} \in \mathcal{A} \\ |A_{1}| = n_{1} \\ m(A_{1}) \geqslant i_{1}}}f_{m(A_{1})-1, i_{1}}f_{A_{1}}\phi_{A_{1}}\Bigg) \dots \Bigg(\sum_{\substack{A_{k} \in \mathcal{A} \\ |A_{k}| = n_{k} \\ m(A_{k}) \geqslant i_{k}}}f_{m(A_{k})-1, i_{k}}f_{A_{k}}\phi_{A_{k}}\Bigg),
\end{align*}
where the second equality follows from the linearity of $L$. Then applying (\ref{recurrence_eq}), extracting coefficients of $b_{i}y^{i}$ from both sides, and rearranging and reindexing, we obtain
\begin{align}
\label{gen_comb_identity_eq}
&\sum_{j=i-1}^{n}f_{j,i}\sum_{\substack{A \in \mathcal{A} \\ |A| = n+1 \\ m(A) = j+1}}f_{A}\phi_{A} - a_{0}(f_{0,0} + f_{0,1})\mathbb{1}_{n=0} \\
\nonumber &\qquad\qquad = \sum_{j=\max\{i-1,0\}}^{n}f_{j,i}\sum_{k=1}^{n}\phi_{k}\sum_{\substack{n_{1} + \cdots + n_{k} = n \\ n_{\ell} \geqslant 1}}\sum_{\substack{i_{1} + \cdots + i_{k} = j \\ 0 \leqslant i_{\ell} \leqslant n_{\ell}}}a_{j}b_{i_{1}}\cdots b_{i_{k}} \\ 
\nonumber &\qquad\qquad\qquad\qquad \times \Bigg(\sum_{\substack{A_{1} \in \mathcal{A} \\ |A_{1}| = n_{1} \\ m(A_{1}) \geqslant i_{1}}}f_{m(A_{1})-1, i_{1}}f_{A_{1}}\phi_{A_{1}}\Bigg) \dots \Bigg(\sum_{\substack{A_{k} \in \mathcal{A} \\ |A_{k}| = n_{k} \\ m(A_{k}) \geqslant i_{k}}}f_{m(A_{k})-1, i_{k}}f_{A_{k}}\phi_{A_{k}}\Bigg)
\end{align}
for all $0 \leqslant i \leqslant n + 1$. This implies recurrence (\ref{recurrence_eq})---if we can prove this identity for some class of combinatorial objects $\mathcal{A}$, then it follows inductively that (\ref{gen_solution_eq}) is the unique solution to the recurrence and therefore the corresponding generating function $G(x, y)$ is the unique solution to equation (\ref{gen_func_eq2}).

\section{Chord diagrams and their features}
\label{chord_diagram_sec}

We now turn to describing two classes of combinatorial objects $\mathcal{A}$ which satisfy an identity of the form (\ref{gen_comb_identity_eq}), thereby formulating our generating function solutions to equations (\ref{bin_func_eq}) and (\ref{power_func_eq}). This requires defining the relevant combinatorial objects, connected chord diagrams, as well as various properties and features of these objects, including some which will not be needed until later in the paper. 

\begin{definition}[Chord diagram, intersection graph, crossing, nesting]
A {\em rooted chord diagram $C$} of size $n$ is a perfect matching of $\{1, 2, \ldots, 2n\}$; the chord containing $1$ is the {\em root} of $C$ and we view the empty set $\emptyset$ as a diagram of size 0. By convention, the elements of $C$ are ordered pairs $(x, y)$ with $x < y$; $x$ and $y$ are called the {\em source} and {\em sink}, respectively, of chord $(x, y)$. A subset $D \subseteq C$ is a {\em subdiagram} of $C$. The {\em directed intersection graph $G(C)$ of $C$} has the chords as vertices and two chords $c_{1} = (x_{1}, y_{1})$ and $c_{2} = (x_{2}, y_{2})$ with $x_{1} < x_{2}$ joined by a directed edge $c_{1}c_{2}$ if $x_{2} < y_{1} < y_{2}$. In this case, we say that $c_{2}$ is a {\em right neighbor of $c_{1}$} and $c_{1}$ is a {\em left neighbor of $c_{2}$}. Forgetting direction, the edge $c_{1}c_{2}$ is also referred to as a {\em crossing} and $c_{1}$ and $c_{2}$ are said to {\em cross}. On the other hand, if $x_{1} < x_{2} < y_{2} < y_{1}$ then $c_{2}$ is {\em nested under $c_{1}$} and together these two chords form a {\em nesting}.  
\end{definition}

Note that we treat a subdiagram $D$ of $C$ as a proper chord diagram on $\{1, 2, \ldots, 2|D|\}$ by applying an order-preserving bijection. In the sequel, any graph-theoretic notions used in the context of a chord diagram or its elements should be understood as referring to the appropriate feature of its directed intersection graph (e.g. $c_{1}$ is in the neighborhood of $c_{2}$). Furthermore, for brevity, we will use the generic term `diagram' to refer to a chord diagram. We will also extend the notion of nesting to subdiagrams in the obvious way: $D' \subseteq C$ is nested under $D \subseteq C$ if each chord of $D'$ is nested under every chord of $D$. Unsurprisingly, a diagram is called {\em nonnesting} (resp. {\em noncrossing}) if it contains no nestings (resp. crossings).  

\begin{figure}[t]
\begin{center}
\includegraphics[scale=0.8]{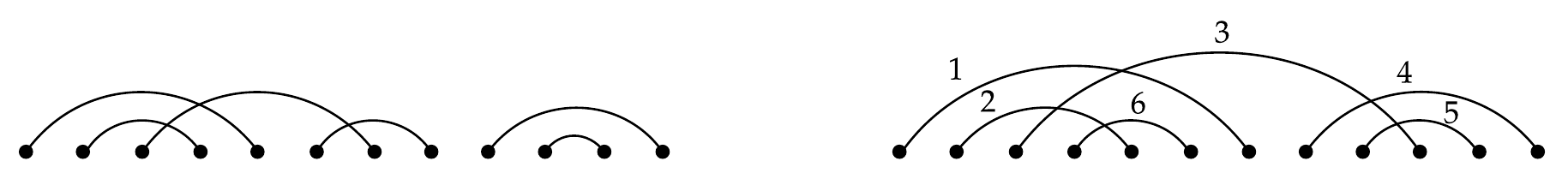}
\caption{Left: a linear representation of a decomposable chord diagram with two indecomposable components and three connected components. Right: a connected chord diagram $C$ with three terminal chords with indices 4, 5, and 6 in the intersection order, which differs from the standard order on $C$.}
\label{chord_diagram_fig}
\end{center}
\end{figure}

There are two representations of a chord diagram used in the literature. The circular representation involves arranging points labeled by $1, 2, \ldots, 2n$ on a circle and, for each chord, joining the source and sink points with a straight line. We can obtain the linear representation from the circular one by cutting the circle just before the root point labeled 1 and straightening the resulting curve, letting the lines corresponding to each chord bend into smooth curves lying above the straightened curve, and then deleting the curve (see Figure~\ref{chord_diagram_fig}). While the circulation representation becomes unrooted after dropping the labels, the linear representation is naturally rooted at the leftmost point and the element labels are determined by the linear order of the points. We will use the linear representation throughout this paper. 

\begin{definition}[Connectivity, connected diagram] The {\em vertex connectivity} $\kappa(C)$, or simply {\em connectivity}, of a diagram $C$ is the vertex connectivity $\kappa(G(C))$ of its intersection graph, that is, the minimum number of vertices whose removal disconnects $G(C)$. The {\em edge connectivity} $\lambda(C)$ of $C$ is defined similarly by removing edges instead of vertices. A diagram $C$ is {\em $k$-connected} if $\kappa(C) \geqslant k$; in particular, it is {\em connected} if it is 1-connected.   
\end{definition}

Note that this definition is really in terms of the undirected intersection graph; the directions on the edges do not play a role. A disconnected diagram can be equivalently defined as one which can be partitioned into two nonempty subdiagrams with no crossings between them. The strictly stronger notion of decomposability of a diagram arises by excluding both crossings and nestings between the two subdiagrams.

\begin{definition}[Diagram concatenation, indecomposable diagram]
Let $C_{1}$ and $C_{2}$ be diagrams. The {\em concatenation} of $C_{1}$ and $C_{2}$ is the diagram $C_{1}C_{2}$ of size $|C_{1}| + |C_{2}|$ whose perfect matching restricts to $C_{1}$ on the first $2|C_{1}|$ elements and to $C_{2}$ on the next $2|C_{2}|$ elements (as subdiagrams). A diagram $C$ is {\em decomposable} if it can be expressed as the concatenation of two smaller diagrams and {\em indecomposable} otherwise. 
\end{definition}

Figure~\ref{chord_diagram_fig} illustrates examples of a decomposable diagram and a connected diagram. Analogously to connected components, we refer to maximal nonempty indecomposable subdiagrams as indecomposable components. Observe the following basic fact that we will use heavily for proofs involving induction. 

\begin{lemma}
The diagram obtained by removing the root chord of a connected diagram is indecomposable. 
\end{lemma}

While the above definitions are largely standard in the combinatorics literature, the following are more unique to this context and most of them first appeared in a paper of Marie and Yeats \cite{Marie2013} (see also \cite{Hihn2019, Courtiel2017a, Courtiel2017b, Courtiel2019}). Note that the directed intersection graph is acyclic, so it induces a partial order on the chords by reachability. We define several different total orders on the chords of a diagram $C$ which extend this partial order. 

\begin{definition}[Standard order, intersection order]
The {\em standard order} of $C$ is given by the order of the sources of the chords of $C$. For connected $C$, the {\em intersection order} is defined recursively as follows: starting with 1, label the root chord of $C$ with the next available label, then remove the root and label the resulting connected components recursively by the standard order of their roots.  
\end{definition}

The standard order will be used as the default order on a chord diagram unless stated otherwise. Note that the standard order and intersection order generally differ substantially (see Figure~\ref{chord_diagram_fig}). If we replace `connected' with `indecomposable' and label the components by the reverse standard order of their roots then we get the {\em peeling order}, which for some purposes is in fact equivalent to the intersection order (see \cite{Courtiel2017b} for details). The weighted generating functions solving (\ref{bin_func_eq}) and (\ref{power_func_eq}) will principally depend on the the positions or indices of certain special chords in the intersection order.  

\begin{definition}[$k$-terminal chords and diagrams]
Let $C$ be a chord diagram. For $k \in \mathbb{N}$, a chord $c \in C$ is {\em $k$-terminal} if it is incident to at most $k-1$ outgoing edges in $G(C)$, that is, it has at most $k-1$ right neighbors.  We refer to 1-terminal chords as simply {\em terminal}; these are the maximal elements in the reachability poset. We extend this language to diagrams: $C$ is {\em $k$-terminal} if there is no $j$-terminal chord before the $j$th-to-last chord for all $1 \leqslant j \leqslant k$. Furthermore, the {\em terminality of $C$} is $k$ if $C$ is $k$-terminal but not $(k+1)$-terminal. 
\end{definition}

Since the chord with rightmost sink in each connected component is necessarily terminal, it follows immediately from the definition that $k$-terminal diagrams are connected and the last $k+1$ chords form a clique; in particular, there are exactly $j$ $j$-terminal chords for all $1 \leqslant j \leqslant k$. Furthermore, clearly $k$-terminal chords and diagrams are also $j$-terminal for $1 \leqslant j \leqslant k$. The indices of the terminal chords will play an important role in both our generating function solutions and the rest of the paper. Accordingly, we will denote the index of the $j$th terminal chord in the intersection order of a diagram $C$ by $t_{j}(C)$. We now record a series of basic facts about these orders and 1-terminal objects, most of which do not seem to have appeared in the existing literature. For the rest of this section let $C$ be a connected diagram of size $n$ with $t_{1}(C) = k$ and $c_{1} < c_{2} < \cdots < c_{n}$ be the chords of $C$ in the intersection order. Note that if $t_{1}(C) = 1$, that is, the root chord is terminal, then it is the only chord of $C$, while if $t_{1}(C) = n$ then $C$ is 1-terminal. 

\begin{lemma}
\label{order_agree_lem}
In the standard order, we also have $c_{1} < c_{2} < \cdots < c_{k}$; that is, the intersection order and standard order agree in a relative sense up to the first terminal chord, which is the chord with rightmost sink. 
\end{lemma}
\begin{proof}
Note that the root chord comes first in the standard order and if we remove it to obtain an indecomposable diagram $C'$, $c_{k}$ is also the first terminal chord in the outermost component of $C'$. With these facts the result follows by a simple induction. 
\end{proof}

The second part of this lemma was observed by Courtiel and Yeats \cite{Courtiel2017a} in their work on terminal chords. As a consequence, the intersection order and standard order are equivalent on 1-terminal diagrams. The other direction does not hold; there are diagrams with multiple terminal chords in which the two orders agree (e.g. take a root crossing a set of terminal chords which all pairwise nest). 


\begin{lemma}
\label{comp_nbhr_lem}
The indecomposable components $C_{1}, C_{2}, \ldots, C_{m}$ remaining after removing $c_{1}, c_{2},\allowbreak \ldots, c_{k}$ have no right neighbors in $C$. 
\end{lemma}
\begin{proof}
Since the intersection order extends the partial order on chords induced by the directed intersection graph, the chords $c_{1}, \ldots, c_{k}$ can only cross chords in $C_{1} \cup \cdots \cup C_{m}$ on the left. 
\end{proof}

\begin{figure}[t]
\begin{center}
\includegraphics[scale=0.8]{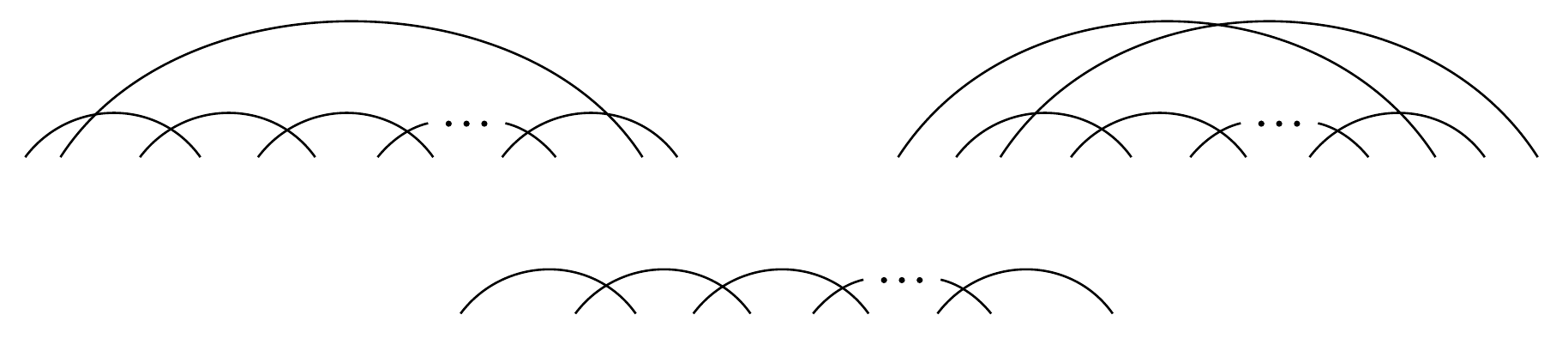}
\caption{Top: top and bottom cycles, the only diagrams whose intersection graph form an induced cycle. Bottom: the unique diagram whose intersection graph forms a nonnesting induced path.}
\label{unique_represent_fig}
\end{center}
\end{figure}

Note that there are exactly two diagrams whose undirected intersection graph is isomorphic to an induced cycle of size $n$; we call these the {\em top cycle} diagram and the {\em bottom cycle} diagram. This diagram representation (near-)uniqueness property was first observed by Bouchet \cite{Bouchet1994}. While no such property holds for induced paths, we gain representation uniqueness if we require that the path is also nonnesting (see Figure~\ref{unique_represent_fig}); this leads to the following observation. 


\begin{lemma}
\label{1term_subdiagram_lem}
There exists a nonnesting induced path in $\{c_{1}, c_{2}, \ldots, c_{k}\}$ from $c_{j}$ to $c_{k}$ for all $1 \leqslant j \leqslant k$. In particular, the chords $c_{1}, c_{2}, \ldots, c_{k}$ induce a 1-terminal subdiagram. 
\end{lemma}
\begin{proof}
If $n = 1$, the result holds trivially. Otherwise, remove the root and consider the outermost component $D$ of the resulting indecomposable diagram. Inductively, for all $2 \leqslant j \leqslant k$ there is a nonnesting induced path $P_{j} \subseteq \{c_{2}, \ldots, c_{k}\}$ from $c_{j}$ to $c_{k}$ in $D$. Since $C$ is connected, by the order extension property the source of $c_{2}$ lies before the sink of $c_{1}$. Furthermore, $c_{k}$ either crosses the root or its source lies after the sink of the root. Thus some chord in $P_{2}$ crosses $c_{1}$; choose such a chord $c_{i}$ with $i$ maximum and let $P'_{2}$ be the subpath of $P_{2}$ beginning at $c_{i}$. Then $\{c_{1}\} \cup P'_{2}$ is a nonnesting induced path from $c_{1}$ to $c_{k}$ in $C$. 
\end{proof}

One can check that this proof actually implicitly constructs nonnesting induced paths to the first terminal chord defined by the property that the $(i+1)$th chord of the path is the rightmost neighbor in $\{c_{1}, c_{2}, \ldots, c_{k}\}$ of the $i$th chord. It also gives a characterization of 1-terminality. 

\begin{corollary}
\label{1term_char_cor}
The following are equivalent.
\begin{enumerate}
\item[(i)] $C$ is 1-terminal.
\item[(ii)] $C - c_{1}$ is 1-terminal. 
\item[(iii)] There exists a nonnesting induced path from $c$ to $c_{n}$ for all $c \in C$. 
\end{enumerate}
\end{corollary}
\begin{proof}
Lemma~\ref{1term_subdiagram_lem} gives the equivalence of (i) and (iii), while the fact that the root $c_{1}$ is not a right neighbor of any chord straightforwardly implies that (i) and (ii) are equivalent. 
\end{proof}

\begin{figure}[t]
\begin{center}
\includegraphics[scale=0.75]{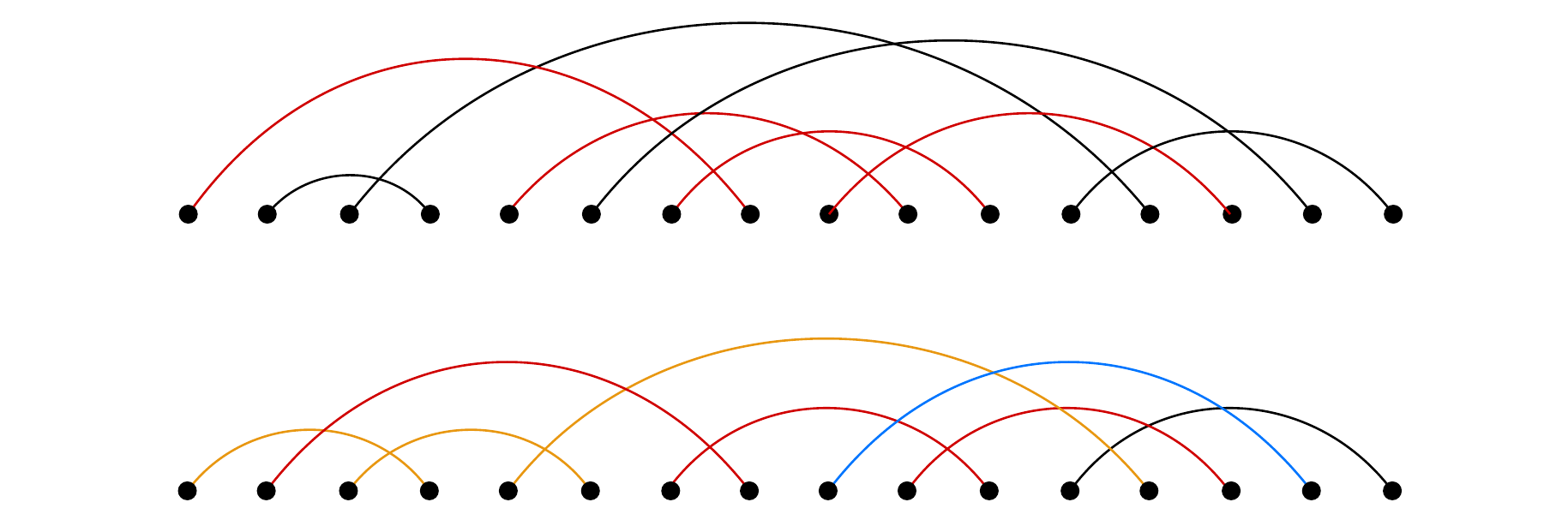}
\caption{Above: a 1-terminal diagram with the traced subdiagram of the rightmost red chord indicated in red. Below: another 1-terminal diagram with the traced subdiagrams of the neighbors of the terminal chord colored. There is no nonnesting induced path between the second orange chord and the second red chord.}
\label{1term_traced_fig}
\end{center}
\end{figure}

The above results begin to illustrate how 1-terminality can be thought of as a more well-behaved notion of connectivity. In particular, for vertex connectivity, a version of (ii) no longer applies and we drop ``nonnesting" to get a version of (iii). We will see further evidence for this view later in the paper. Figure~\ref{1term_traced_fig} illustrates a 1-terminal diagram with two chords that are not linked by a nonnesting induced path, indicating that this result cannot be strengthened to get such a path between e.g. any two nonnesting chords absent a stronger hypothesis. 

Our proofs in the next section involve decomposing a connected chord diagram into specific subdiagrams. This decomposition is constructed by first decomposing the 1-terminal part of the diagram revealed by Lemma~\ref{1term_subdiagram_lem} into certain 1-terminal subdiagrams, which we define presently. We begin with a key concept useful to defining these objects as well as further aspects of decompositions and maps introduced later in the paper. 

\begin{definition}[Source-sink group]
\label{source_sink_group_def}
We say that a chord $c \in C$ is {\em attached} to an indecomposable subdiagram $B$ of $C$ if $\{c\}$ is not a connected component of $B \cup \{c\}$; symmetrically, $B$ is also {\em attached} to $c$. Define the {\em source-sink group} of $c$ to be the interval of endpoints (viewed as a part of the diagram) containing the source of $c$, the maximal interval of sinks of $\{c_{1}, \ldots, c_{k}\}$ immediately following it excluding the sink of $c$, and the indecomposable components of $C - \{c_{1}, \ldots, c_{k}\}$ attached to any chord containing one of these sinks. 
\end{definition}

More informally, we will also refer to chords as ``attached" to their sources and sinks. Let $D$ be a 1-terminal diagram. By Lemma~\ref{comp_nbhr_lem} the source-sink groups of $D$ contain only the source and adjacent maximal interval of sinks. 

\begin{definition}[Traced subdiagram]
\label{traced_subdiagram_def}
For $c \in D$, we construct the {\em traced subdiagram $D_{c}$ of $c$} as follows. The chord $c$ belongs to $D_{c}$. Beginning with $c' = c$, place in $D_{c}$ the chord attached to every sink in the source-sink group of $c'$, not including the sink of $c'$, and then repeat this procedure with each newly added chord of $D_{c}$.  
\end{definition}

See Figure~\ref{1term_traced_fig} for examples. Note that there may be no sinks in the source-sink group of $c'$; in particular, by connectivity, this holds in the case in which either $c'$ is the only chord of $D$ or there is another source next to the source of $c'$. Since every chord attached to a sink in the source-sink group crosses the chord of the source to the left, traced subdiagrams are equivalently specified as follows: each chord of $D$ besides the base chord $c$ is contained in $D_{c}$ if and only if its last right neighbor (in the standard order) is in $D_{c}$. In particular, we observe the following. 

\begin{lemma}
\label{traced_1term_lem}
Traced subdiagrams $D_{c}$ are 1-terminal and $c$ is the terminal chord. 
\end{lemma}

The above observation together with the following result gives a unique decomposition of a 1-terminal diagram into smaller 1-terminal subdiagrams, after removing the terminal chord. 

\begin{lemma}
\label{traced_partition_lem}
The traced subdiagrams of the neighbors of the terminal chord $d$ of $D$ partition $D - d$. 
\end{lemma}
\begin{proof}
If $|D| = 1$, the result trivially holds. Otherwise, the root chord $c$ is either a neighbor of $d$ or it is not. If it is, note that by 1-terminality there are no chords nested under the terminal chord. Then clearly $D_{c} = \{c\}$ and, since the last right neighbor of $c$ is $d$, it is not contained in any other traced subdiagram besides $D_{d}$. Since traced subdiagrams are 1-terminal, by Corollary~\ref{1term_char_cor} it follows that all other traced subdiagrams (excluding $D_{d}$) are inherited from $D - c$. Then we get the required partition inductively or immediately if $D - c - d$ is empty. On the other hand, if $c$ is not a neighbor of $d$ then it is contained in the traced subdiagram of $D$ of its rightmost neighbor, which exists since $D$ is connected. Then, as before, inductively we get the required partition. 
\end{proof}


Although we will not need this fact, it is also worth observing that this implies that the traced subdiagram of the terminal chord is the entire diagram.

\subsection{Weighting connected chord diagrams}

We require one final piece of the puzzle: weights for connected chord diagrams. Analogously to trees, we begin by defining a notion of {\em valency} for chords in a diagram. Let $C$ be a connected chord diagram. 

\begin{definition}
\label{valency_def}
Define the {\em valency} of $c \in C$ to be $k + \ell$ if there are exactly $k$ left neighbors of $c$ not attached to any chord $c' > c$ and after removing the left neighbors of $c$ we are left with exactly $\ell$ indecomposable components immediately proceeding the source of $c$. 
\end{definition}

\begin{figure}[t]
\begin{center}
\includegraphics[scale=0.85]{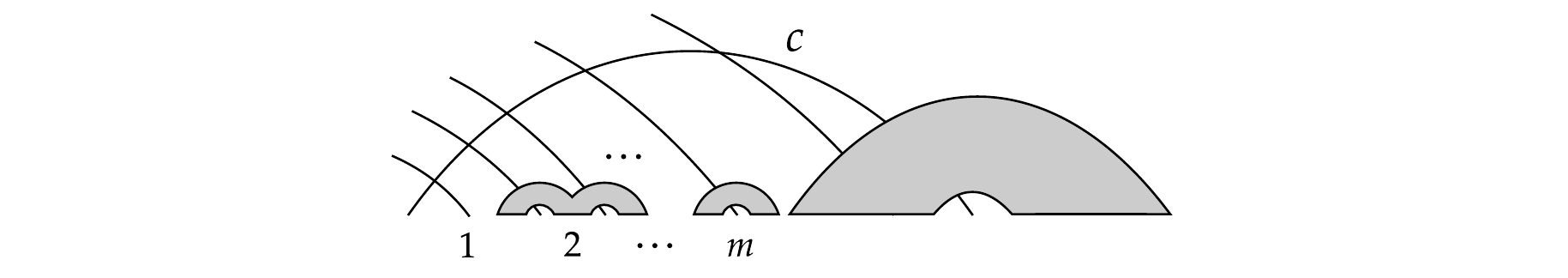}
\caption{A chord $c$ of valency $m$ in a connected diagram.}
\label{chord_valency_fig}
\end{center}
\end{figure}

Loosely speaking, in a local sense the valency is the number of connected pieces in the source-sink group of a chord. This notion of valency for a diagram should be thought of as playing the same role as outdegree for a tree, that is, the number of children of a vertex in a rooted tree. Figure~\ref{chord_valency_fig} provides a generic visual representation of a chord of valency $m$. We can easily see that if $C$ is 1-terminal then the valency of $c' \in C$ is exactly the number of chords added at the corresponding step in the construction of a traced subdiagram containing $c'$. 

Given a sequence of weights $(\phi_{k})_{k \geqslant 0}$ in a field $K$, we associate a weight $\phi_{k}$ with each chord $c \in C$ of valency $k$. Then we define the {\em weight} of $C$ to be 
\begin{align*}
\phi_{C} = \prod_{c \in C}\phi_{\textrm{val}(c)}.
\end{align*}
This is exactly the same way weights of trees are traditionally defined, at least in the context of tree models and hook lengths \cite{Flajolet2009, Kuba2013}. We will elaborate more on the relationship between weighted trees and weighted chord diagrams later in the paper. 

For now, we briefly describe how this relates to the weights which arise in the more physically-grounded work of Hihn and Yeats \cite{Hihn2019} and Courtiel, Yeats, and Zeilberger \cite{Courtiel2017b} on generalized Dyson-Schwinger equations. In the latter, these weights are cast as a product of binomial coefficients ${\omega_{i}(C) + s - 1 \choose \omega_{i}(C)}$ over all chords $c_{i}$ in a diagram $C$, where $\omega_{i}(C)$ is a parameter known as the covering number of $c_{i}$ and $s$ is a positive integer. The covering number of $c_{i}$ is defined as the number of intervals, that is, gaps between points of the diagram, assigned to $c_{i}$ in a certain recursive covering procedure. One can readily see that these intervals correspond exactly to the intervals on either side of each of the $m$ attachment points and components from Definition~\ref{valency_def} of a chord of valency $m$; in other words, $ \textrm{val}(c) = \omega_{i}(C) - 1$. So setting $\phi_{k} = {k + s \choose k + 1}$ gives the weighted cases considered in prior work. Note though that \cite{Hihn2019} and \cite{Courtiel2017b} also simultaneously apply a separate weighting to the chords, so our work is not a full generalization of these studies.

\section{The generating function solutions}
\label{func_eq_sol_sec}

We are now ready to formulate and, in the next two subsections, prove our main result. 

\begin{theorem}
\label{equation_sol_thm}
The functional equation 
\begin{align}
\label{general_func_eq2}
G(x, y) = xL(\phi(G(x, y)))
\end{align}
is uniquely solved by 
\begin{align}
\label{binom_sol_eq}
G(x, y) = \sum_{C \in \mathcal{C}}f_{C}\phi_{C}x^{|C|}{L_{bin}(y^{t_{1}(C)-1}) \over (t_{1}(C)-1)!} \quad \textrm{if} \quad L = L_{bin}
\end{align}
and 
\begin{align}
\label{power_sol_eq}
G(x, y) = \sum_{C \in \mathcal{C}_{top}}f_{C}\phi_{C}x^{|C|}L_{div}(y^{t_{1}(C)-1}) \quad \textrm{if} \quad L = L_{div},
\end{align}
where $\mathcal{C}$ is the set of connected chord diagrams, $\mathcal{C}_{top}$ is the set of connected top-cycle-free chord diagrams, that is, connected diagrams with no top cycle subdiagram, and 
\begin{align*}
f_{C} = f_{0}^{|C| - k}f_{t_{2}(C)-t_{1}(C)}f_{t_{3}(C)-t_{2}(C)} \cdots f_{t_{k}(C)-t_{k-1}(C)}
\end{align*}
for a diagram $C$ with $k$ terminal chords and $f_{i}$ the coefficients of the power series $F$ specifying $L$. 
\end{theorem}

In other words, when $L$ is a polynomial 1-cocycle operator (\ref{general_func_eq2}) is solved by a generating function for a certain family of connected weighted chord diagrams  where each term is weighted by a monomial in the coefficients defining $L$ determined by the size of the diagram, the number of terminal chords, and the differences between the indices of consecutive terminal chords in the intersection order. The diagrams are counted by their size in the $x$ variable and one less than the index of the first terminal chord in the $y$ variable. While (\ref{power_sol_eq}) can be thought of as ordinary in both variables, (\ref{binom_sol_eq}) should be viewed as ordinary in $x$ and exponential in $y$, although the fact that $L(y^{t_{1}(C)-1})$ appears in the series rather $y^{t_{1}(C)-1}$ makes this less than strictly true. One could account for this by instead regarding both series as ordinary generating functions for a set of polynomials $\{p_{C}(y)\}$ indexed by a class of weighted chord diagrams. It is notable that the monomial $f_{C}$ is determined by the gaps between terminal chord indices, not the indices themselves, indicating that it is their relative position that matters, not their absolute position.

\subsection{The binomial 1-cocycle equation}
\label{bin_eq_sec}

We begin by fixing $L = L_{bin}$ and proving that (\ref{binom_sol_eq}) uniquely solves (\ref{general_func_eq2}). 
Note that applying $L$ to the standard basis of $K[y]$ gives
\begin{align*}
L(y^{n}) = \int_{0}^{y}F{\left({d \over dt}\right)}t^{n}dt &= \int_{0}^{y}\sum_{i \geqslant 0}f_{i}{d^{i} \over dt^{i}}t^{n}dt \\
&= \int_{0}^{y}\sum_{i=0}^{n}f_{i}{n! \over (n-i)!}t^{n-i}dt \\
&= n!\sum_{i=1}^{n+1}f_{n+1-i}{y^{i} \over i!}.
\end{align*}
Then, letting $\mathcal{A}$ be the set $\mathcal{C}$ of weighted connected chord diagrams, $m(C)$ be the index $t_{1}(C)$ of the first terminal chord of $C \in \mathcal{C}$, $a_{n} = n!$, $b_{i} = {1 \over i!}$, and $f_{n,i} = f_{n+1-i}$ if $i > 0$ and 0 otherwise, from (\ref{gen_comb_identity_eq}) it suffices to show that 
\begin{align*}
\sum_{j=i-1}^{n}f_{j+1-i}\sum_{\substack{C \in \mathcal{C} \\ |C| = n+1 \\ t_{1}(C) = j+1}}f_{C}\phi_{C} - f_{0}\mathbb{1}_{n=0} &= \sum_{j=\max\{i-1, 1\}}^{n}f_{j+1-i}\sum_{k=1}^{n}\phi_{k}\sum_{\substack{n_{1} + \cdots + n_{k} = n \\ n_{\ell} \geqslant 1}}\sum_{\substack{i_{1} + \cdots + i_{k} = j \\ 1 \leqslant i_{\ell} \leqslant n_{\ell}}}{j \choose i_{1}, \cdots, i_{k}} \\
&\quad \times \Bigg(\sum_{\substack{C_{1} \in \mathcal{C} \\ |C_{1}| = n_{1} \\ t_{1}(C_{1}) \geqslant i_{1}}}f_{t_{1}(C_{1})-i_{1}}f_{C_{1}}\phi_{C_{1}}\Bigg) \cdots \Bigg(\sum_{\substack{C_{k} \in \mathcal{C} \\ |C_{k}| = n_{k} \\ t_{1}(C_{k}) \geqslant i_{k}}}f_{t_{1}(C_{k})-i_{k}}f_{C_{k}}\phi_{C_{k}}\Bigg)
\end{align*}
for all $1 \leqslant i \leqslant n+1$. Observe that the difference in the left hand side series with $i = 1$ and $i = 2$ is the contribution from the connected diagrams in which the first chord is terminal, of which there is only the single diagram of size one. Since it is the unique diagram with $f_{C} = 1$ and $\phi_{C} = 1$, the corresponding term in the series cancels with the term $f_{0}\mathbb{1}_{n=0}$ on the left hand side. It follows then that we may reindex and thereby remove the sums over $j$ on both sides to obtain
\begin{align}
\label{bin_identity_eq}
\sum_{\substack{C \in \mathcal{C} \\ |C| = n+1 \\ t_{1}(C) = j+1}}f_{C}\phi_{C}&= \sum_{k=1}^{n}\phi_{k}\sum_{\substack{n_{1} + \cdots + n_{k} = n \\ n_{\ell} \geqslant 1}}\sum_{\substack{i_{1} + \cdots + i_{k} = j \\ 1 \leqslant i_{\ell} \leqslant n_{\ell}}}{j \choose i_{1}, \cdots, i_{k}} \\
\nonumber &\qquad\qquad \times \Bigg(\sum_{\substack{C_{1} \in \mathcal{C} \\ |C_{1}| = n_{1} \\ t_{1}(C_{1}) \geqslant i_{1}}}f_{t_{1}(C_{1})-i_{1}}f_{C_{1}}\phi_{C_{1}}\Bigg) \cdots \Bigg(\sum_{\substack{C_{k} \in \mathcal{C} \\ |C_{k}| = n_{k} \\ t_{1}(C_{k}) \geqslant i_{k}}}f_{t_{1}(C_{k})-i_{k}}f_{C_{k}}\phi_{C_{k}}\Bigg).
\end{align}
It thereby remains to show that this identity holds for $1 \leqslant j \leqslant n$. To prove this identity, we will show how it is a statement about the existence of a certain unique decomposition of a connected chord diagram into several smaller shuffled connected diagrams. More formally, it suffices to construct a map $\alpha$ from connected diagrams $C$ with $t_{1}(C) = j + 1$ to $k$-tuples of diagrams satisfying the required constraints on size and index of the first terminal chord with each diagram associated to a non-empty subset in a partition of $[1, j]$, exhibit the inverse of $\alpha$, and check that the weights and monomials of $f_{i}$'s behave correctly under $\alpha$. Since the number of partitions of $[1, j]$ with blocks of size $i_{1}, \ldots, i_{k}$ is given by the multinomial coefficient ${j \choose i_{1}, \ldots, i_{k}}$, (\ref{bin_identity_eq}) will follow. 

So, let $C$ be a connected diagram of size $n+1$ with $t_{1}(C) = j + 1$ and let $c_{1} < c_{2} < \cdots < c_{n+1}$ be the chords of $C$ in the intersection order. Then we define the following. 
\begin{itemize}
\item $D$: the 1-terminal diagram induced by chords $c_{1}, \ldots, c_{j+1}$ of $C$.
\item $A_{1}, A_{2}, \ldots, A_{k}$: the indecomposable components of $C - \{c_{1}, \ldots, c_{j}\}$ nested under $c_{j+1}$. Neighbors of the terminal chord $c_{j+1}$ may or may not be attached to such a component; for each that is not, we also include among these $A_{\ell}$'s an empty component. 
\item For all $1 \leqslant \ell \leqslant k$, $D_{\ell}$: the union of the traced subdiagrams of $D$ of each neighbor of $c_{j+1}$ attached to component $A_{\ell}$ in $C$; if $A_{\ell}$ is empty this is simply the traced subdiagram of the neighbor defining $A_{\ell}$ above. 
\item For all $1 \leqslant \ell \leqslant k$, $C_{\ell}$: the union of $D_{\ell}$ with all indecomposable components of $C - D$ attached to chords of $D_{\ell}$ in $C$.
\item For all $1 \leqslant \ell \leqslant k$, $I_{\ell}$: the subset $\{1 \leqslant i \leqslant j \mid c_{i} \in C_{\ell}\}$ of $[1, j]$ indicating which chords in $D - c_{j+1}$ belong to $C_{\ell}$. Furthermore, set $i_{\ell} = |I_{\ell}|$.
\item $\alpha(C) = ((C_{1}, I_{1}), (C_{2}, I_{2}), \ldots, (C_{k}, I_{k}))$. 
\end{itemize}

\begin{figure}[t]
\begin{center}
\includegraphics[scale=0.8]{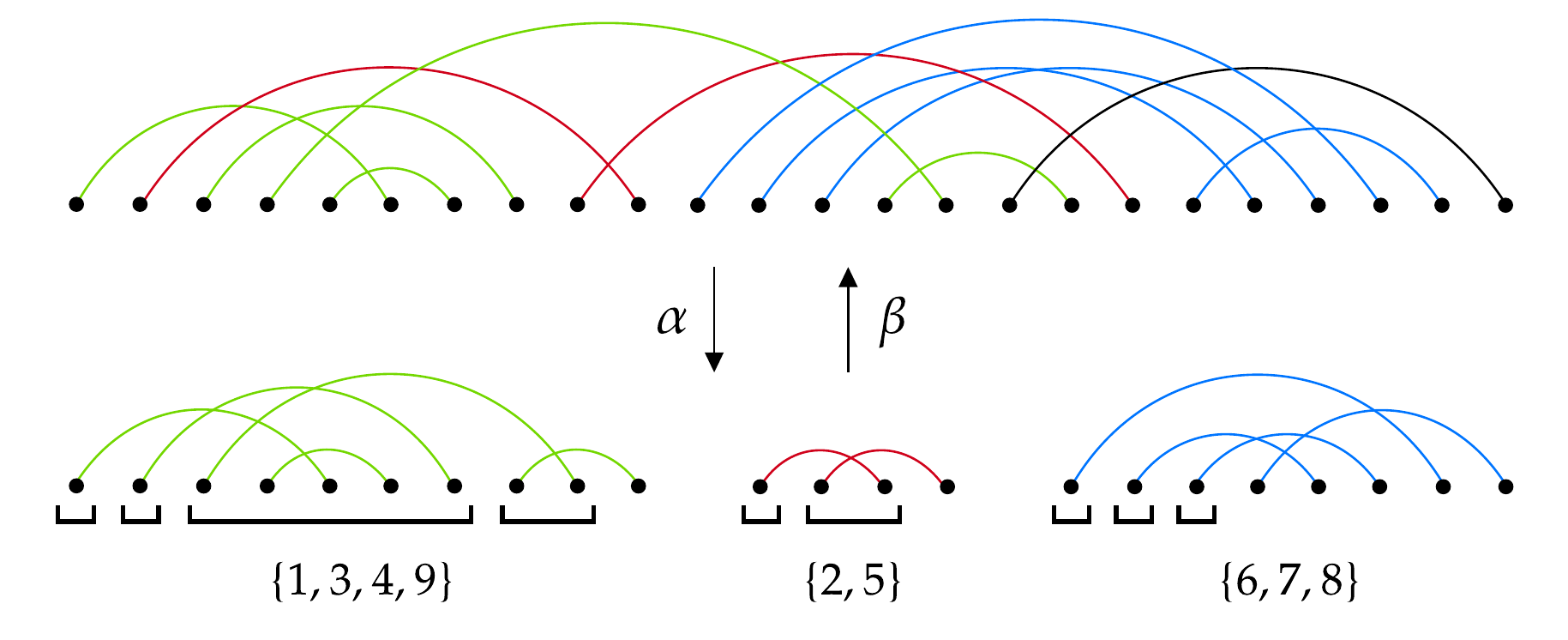}
\caption{A connected diagram and its permuted decomposition defining the maps $\alpha$ and $\beta$. The source-sink groups used in the construction of $\beta$ are indicated by the brackets below the three diagrams on the right side of the figure.}
\label{alpha_beta_decomp_fig}
\end{center}
\end{figure}

Figure~\ref{alpha_beta_decomp_fig} illustrates an example of this construction for a representative diagram $C$. By Lemma~\ref{comp_nbhr_lem} chords of $D$ attach to indecomposable components of $C - D$ at sinks; in other words, the interval occupied by an indecomposable component of $C - D$ intersects at most one maximal interval of sinks of $D$. By the construction of traced subdiagrams, it follows that each indecomposable component of $C - D$ is attached to the chords of $D_{\ell}$ for at most one $\ell$. Thus, we observe that $\{C_{1}, \ldots, C_{k}\}$ is a partition of $C - c_{j+1}$ by Lemma~\ref{traced_partition_lem} and the fact that $C$ is connected (since therefore every indecomposable component of $C - D$ is attached to a chord of $D - c_{j+1}$); note that this also implicitly uses the fact that every neighbor of $c_{j+1}$ is in $D$ since the intersection order is an extension of the reachability order. We also have the following properties. 

\medskip

\noindent
\textbf{Claim.} For all $\ell$, (i) $C_{\ell}$ is connected and (ii) $t_{1}(C_{\ell}) \geqslant i_{\ell}$. 

\noindent
\textit{Proof.} Note that clearly each neighbor of $c_{j+1}$ attached to $A_{\ell}$ has to every other attached neighbor since each cross a chord in the outermost connected component of $A_{\ell}$ by definition. Furthermore, the traced subdiagram of $D$ of a neighbor of $c_{j+1}$ is connected since it is 1-terminal by Lemma~\ref{traced_1term_lem}. By construction, it follows that $C_{\ell}$ is connected since $C$ is connected, giving (i). For (ii), note that either $C_{\ell}$ is constructed from the traced subdiagram of exactly one neighbor $c_{i}$ of $c_{j+1}$ or $A_{\ell}$ is nonempty. In the former case, $c_{i}$ is the first terminal of $C_{\ell}$, while in the latter case each chord of $D_{\ell}$ has a right neighbor in $C_{\ell}$ by the 1-terminality of traced subdiagrams. Thus, in either case, each neighbor of $c_{j+1}$ in $C_{\ell}$ precedes (non-strictly) the first terminal chord of $C_{\ell}$, which is necessarily either the neighbor $c_{i}$ or in $A_{\ell}$. Then another application of Lemma~\ref{traced_1term_lem} yields the inequality $t_{1}(C_{\ell}) \geqslant i_{\ell}$.  \qed

\medskip

Clearly the valency of the terminal chord $c_{j+1}$ is $k$. Furthermore, for each chord $c \in C_{\ell}$ the valency of $c$ in $C$ and $C_{\ell}$ is equal by contruction. Then it follows that $\phi_{C} = \phi_{k}\phi_{C_{1}} \cdots \phi_{C_{k}}$. 

Before constructing the inverse of $\alpha$, it remains only to ensure equality of the monomials corresponding to $C$ and $(C_{1}, \ldots, C_{n})$. Let $B_{1}, \ldots, B_{m}$ be the components of $C - D$ listed in the intersection order of $C$; that is, all chords of $B_{1}$ come before all chords of $B_{2}$, and so on. We require a lemma which says that $f_{C}$ is determined by $f_{D}, f_{B_{1}}, \ldots, f_{B_{m}}$ and the index of the first terminal chords of $B_{1}, \ldots, B_{m}$ in the intersection order of $C$. 

\begin{lemma}
\label{monom_factor_lem}
We have
\begin{align*}
f_{C} = f_{D}f_{t_{1}(B_{1})}f_{B_{1}} \cdots f_{t_{1}(B_{m})}f_{B_{m}}.
\end{align*}
In other words, $f_{C} = f_{C - B_{m}}f_{t_{1}(B_{m})}f_{B_{m}}$. 
\end{lemma}
\begin{proof}
By Lemma~\ref{comp_nbhr_lem}, terminal chords of $B_{i}$ are terminal in $C$ and vice versa for all $i$. Furthermore, by Lemma~\ref{1term_subdiagram_lem} $D$ is 1-terminal both as a diagram and as a subdiagram of $C$. Thus every terminal chord of $C$ corresponds bijectively to a terminal chord in $D, B_{1}, \ldots, B_{k}$. Furthermore, clearly the intersection orders on $D, B_{1}, \ldots, B_{k}$ agree with the intersection order on $C$ (in the sense that $c < c'$ in $D$ or $B_{i}$ if and only if $c < c'$ in $C$ for $c, c' \in D$ or $c, c' \in B_{i}$ for some $i$) by construction. It follows that the number of terminal chords of $C$ is equal to the sum of the number of terminal chords in $D$ and $B_{1}, \ldots, B_{k}$. We also further infer by construction that every difference $t_{i+1}(C) - t_{i}(C)$ corresponds uniquely to either a difference of consecutive terminal chord indices in some $B_{i}$ or the difference between the index in $C$ of the last terminal chord of $B_{i-1}$ and the index in $C$ of the first terminal chord of $B_{i}$, that is, $t_{1}(B_{i})$. The desired equality follows from these two observations. 
\end{proof}

This gives the needed equality. 

\begin{lemma}
We have 
\begin{align*}
f_{C} = f_{t_{1}(C_{1})-i_{1}}f_{C_{1}} \cdots f_{t_{1}(C_{k})-i_{k}}f_{C_{k}}.
\end{align*}
\end{lemma}
\begin{proof}
By Lemma~\ref{monom_factor_lem} applied to $f_{C}$ and $f_{C_{1}}, \ldots, f_{C_{n}}$, we may assume that $C - D - \bigcup_{\ell}A_{\ell}$ is empty, since the corresponding terms in the monomials can simply be cancelled from both sides. Then $C_{\ell} = D_{\ell} \cup A_{\ell}$. By our observations in the proof of the claim above, the chords of $D_{\ell}$ are the first $i_{\ell}$ chords of $C_{\ell}$ in the intersection order. This implies that 
\begin{align*}
f_{C_{\ell}} = f_{D_{\ell} \cup A_{\ell}} = \begin{cases}
f_{0}^{i_{\ell}}A_{\ell} &\textrm{if } A_{\ell} \textrm{ nonempty} \\
f_{0}^{i_{\ell}-1} &\textrm{otherwise},
\end{cases}
\end{align*}
$t_{1}(A_{\ell}) = t_{1}(C_{\ell}) - i_{\ell}$ (which is zero if and only if $A_{\ell}$ is empty), and 
\begin{align*}
f_{D} = f_{0}^{j} = f_{0}^{i_{1} + \cdots + i_{k}} = f_{0}^{i_{1}}\cdots f_{0}^{i_{k}}.
\end{align*}
Then by Lemma~\ref{monom_factor_lem}
\begin{align*}
f_{C} &= f_{D}f_{t_{1}(A_{1})}f_{A_{1}} \cdots f_{t_{1}(A_{k})}f_{A_{k}} \\
&= f_{0}^{i_{1}}\cdots f_{0}^{i_{k}}f_{t_{1}(C_{1})-i_{1}}f_{A_{1}}\cdots f_{t_{1}(C_{k})-i_{k}}f_{A_{k}} \\
&= f_{t_{1}(C_{1})-i_{1}}f_{C_{1}} \cdots f_{t_{1}(C_{k})-i_{k}}f_{C_{k}}. \qedhere
\end{align*}
\end{proof}

Now we turn to the final step in the proof of (\ref{bin_identity_eq}), showing that $\alpha$ has an inverse, that is, defining a map $\beta$ such that $\beta \circ \alpha = \id$. Suppose we are given a tuple $((C_{1}, I_{1}), \ldots, (C_{k}, I_{k}))$ with $C_{1}, \ldots, C_{k}$ connected diagrams and $\{I_{1}, \ldots, I_{k}\}$ a partition of $[1, j]$ with $1 \leqslant |I_{\ell}| \leqslant t_{1}(C_{\ell})$ for all $\ell$ and $j = |I_{1}| + \cdots + |I_{k}|$. We define the composed diagram $C$ as follows. For all $\ell$, let $D_{\ell}$ be the 1-terminal diagram induced by the first $|I_{\ell}|$ chords of $C_{\ell}$ in the intersection order. Let $C'$ be the diagram obtained by concatenating $C_{1}, \ldots, C_{k}$ and then arranging the source-sink groups corresponding to the first $|I_{1}|, \ldots, |I_{k}|$ of $C_{1}, \ldots, C_{k}$, respectively, in order according to the permutation induced by $I_{1}, \ldots, I_{k}$, maintaining the standard order of the chords within each diagram $C_{\ell}$. Finally, let $C$ be obtained from $C'$ by adding a chord $c$ with its sink in the rightmost position and source immediately proceeding the source-sink groups corresponding to the first $j$ chords of $C'$ (in the intersection order) and set $\beta((C_{1}, I_{1}), \ldots, (C_{k}, I_{k})) = C$. A representative example of this construction is given in Figure~\ref{alpha_beta_decomp_fig}.

For all $\ell$, clearly at least one source of $C_{\ell}$ lies to the left of the source of $c$ while the rightmost sink of $C_{\ell}$ lies to the right of the source of $c$. Since $C_{\ell}$ is connected, it follows that $c$ has a neighbor from each $C_{\ell}$, implying that $C$ is connected. Then by Lemma~\ref{order_agree_lem}, $c$ is the first terminal chord of $C$ since its sink lies furthest to the right. Note that nonnesting paths of $C_{\ell}$ are preserved in $C$ by the construction process; since the first terminal chord of $C_{\ell}$ clearly either crosses $c$ or is nested under it, this implies that there is a nonnesting path from each of the first $|I_{\ell}|$ of $C_{\ell}$ to $c$, implying that $t_{1}(C) \geqslant |I_{1}| + \cdots + |I_{\ell}| + 1 = j+1$ by Lemma~\ref{1term_subdiagram_lem}. On the other hand, by the extension property of the intersection order, all the chords of $C_{\ell} - D_{\ell}$ come after the first terminal chord of $C$ in the intersection order. Combined with the fact that only $|I_{\ell}|$ chords of $D_{\ell}$ have their sources to the left of the source of $c$, we see that $t_{1}(C) \leqslant j + 1$ by Lemma~\ref{order_agree_lem}. We therefore infer that $t_{1}(C) = j + 1$, as required.

Finally, we conclude the proof by observing that $\beta(\alpha(C)) = C$ by construction and our prior conclusions about $\alpha$ and $\beta$.  



\subsection{The divided power 1-cocycle equation}
\label{pow_eq_sec}

We now fix $L = L_{div}$ and prove that (\ref{power_sol_eq}) uniquely solves (\ref{general_func_eq2}). In addition to its simplicity and clarity, this illustrates one of the primary advantages of the proof stategy used in the previous section over the method used by Marie and Yeats \cite{Marie2013} for solving the functional equation with the binomial 1-cocycle: it easily adapts with few changes to also handle the divided power 1-cocycle. The changes will end up essentially amounting to dropping factorials from all the equalities of the previous section, which translates into top-cycle-free diagrams {\em uniquely} decomposing into a tuple of smaller top-cycle-free diagrams---the need to provide a shuffle-defining partition to accompany the decomposition effectively disappears. 

Note that $L$ is given on the standard basis 
of $K[y]$ by 
\begin{align*}
L(y^{n}) = yF{\left(d_{y}\right)}y^{n} = \sum_{i=1}^{n+1}f_{n+1-i}y^{i}.
\end{align*}
Then, letting $\mathcal{A}$ be the set $\mathcal{C}_{top}$ of weighted connected top-cycle-free diagrams, $m(C) = t_{1}(C)$ for $C \in \mathcal{C}_{top}$, $a_{n} = b_{i} = 1$, and $f_{n,i} = f_{n+1-i}$ if $i > 0$ and 0 otherwise, by (\ref{gen_comb_identity_eq}) it suffices to show that 
\begin{align*}
\sum_{j=i-1}^{n}f_{j+1-i}\sum_{\substack{C \in \mathcal{C}_{top} \\ |C| = n+1 \\ t_{1}(C) = j+1}}f_{C}\phi_{C} - f_{0}\mathbb{1}_{n=0} &= \sum_{j=\max\{i-1, 1\}}^{n}f_{j+1-i}\sum_{k=1}^{n}\phi_{k}\sum_{\substack{n_{1} + \cdots + n_{k} = n \\ n_{\ell} \geqslant 1}}\sum_{\substack{i_{1} + \cdots + i_{k} = j \\ 1 \leqslant i_{\ell} \leqslant n_{\ell}}} \\
&\quad \times \Bigg(\sum_{\substack{C_{1} \in \mathcal{C}_{top} \\ |C_{1}| = n_{1} \\ t_{1}(C_{1}) \geqslant i_{1}}}f_{t_{1}(C_{1})-i_{1}}f_{C_{1}}\phi_{C_{1}}\Bigg) \cdots \Bigg(\sum_{\substack{C_{k} \in \mathcal{C}_{top} \\ |C_{k}| = n_{k} \\ t_{1}(C_{k}) \geqslant i_{k}}}f_{t_{1}(C_{k})-i_{k}}f_{C_{k}}\phi_{C_{k}}\Bigg)
\end{align*}
for all $1 \leqslant i \leqslant n+1$. By the same reasoning as in the binomial case, this reduces to the equivalent identity
\begin{align}
\label{pow_identity_eq}
\sum_{\substack{C \in \mathcal{C}_{top} \\ |C| = n+1 \\ t_{1}(C) = j+1}}f_{C}\phi_{C} &= \sum_{k=0}^{n}\phi_{k}\sum_{\substack{n_{1} + \cdots + n_{k} = n \\ n_{\ell} \geqslant 1}}\sum_{\substack{i_{1} + \cdots + i_{k} = j \\ 1 \leqslant i_{\ell} \leqslant n_{\ell}}} \\
&\qquad\quad \times \Bigg(\sum_{\substack{C_{1} \in \mathcal{C}_{top} \\ |C_{1}| = n_{1} \\ t_{1}(C_{1}) \geqslant i_{1}}}f_{t_{1}(C_{1})-i_{1}}f_{C_{1}}\phi_{C_{1}}\Bigg) \cdots \Bigg(\sum_{\substack{C_{k} \in \mathcal{C}_{top} \\ |C_{k}| = n_{k} \\ t_{1}(C_{k}) \geqslant i_{k}}}f_{t_{1}(C_{k})-i_{k}}f_{C_{k}}\phi_{C_{k}}\Bigg)
\end{align}
for $1 \leqslant j \leqslant n$. We now turn to proving this identity. Let $C$ be a connected top-cycle-free diagram, that is, a connected diagram with no top cycle subdiagram, of size $n+1$ with $t_{1}(C) = j+1$ and let $c_{1} < c_{2} < \cdots < c_{n+1}$ be the chords of $C$ in the intersection order. As previously, let $\alpha(C) = ((C_{1}, I_{1}), \ldots, (C_{k}, I_{k}))$.  Note that $C_{1}, \ldots, C_{k}$ are clearly top-cycle-free since they are subdiagrams of the top-cycle-free diagram $C$. 

\medskip

\noindent
\textbf{Claim.} $I_{1}, \ldots, I_{k}$ are intervals. 

\noindent
\textit{Proof.} Suppose to the contrary that $I_{\ell}$ is not an interval for some $\ell$. Then there exists a block $I_{\ell'}$ such that $I_{\ell}$ and $I_{\ell'}$ either nest or cross, that is, there are $r, s \in I_{\ell}$ and $t \in I_{\ell'}$ such that $r < t < s$. By construction, each subdiagram $D_{\ell}$ of $C_{\ell}$ contains a neighbor of the first terminal chord $c_{j+1}$ of $C$ and each such neighbor lies before the first terminal chord of $C_{\ell}$ in the intersection order. Since $D_{\ell}$ is connected, there exists a path $P$ in $D_{\ell}$ from $c_{r}$ to $c_{s}$. By construction of $D_{\ell'}$ and Lemma~\ref{1term_subdiagram_lem}, there exists a nonnesting path $Q$ in $D_{\ell'}$ from $c_{t}$ to a neighbor of $c_{j+1}$ in $D_{\ell'}$; let $Q' = Q \cup \{c_{j+1}\}$. Since $r < t < s$, by Lemma~\ref{order_agree_lem} $c_{r} < c_{t} < c_{s}$ in the standard order, implying that $P$ and $Q$ cross, that is, there is a chord $c$ in $P$ and a chord $c'$ in $Q$ such that $c$ and $c'$ cross. Let $P'$ be a nonnesting path in $D_{\ell}$ from $c$ to a neighbor of $c_{j+1}$ in $D_{\ell}$, which again exists by construction and Lemma~\ref{1term_subdiagram_lem}. Then $P' \cup Q'$ contains two internally-disjoint nonnesting paths of length at least three with the same first and last chords (in the standard order), so one can easily see that it contains a top cycle as a subdiagram. But then $C$ is not top-cycle-free, a contradiction. \qed

\medskip

Now fix a tuple $((C_{1}, I_{1}), \ldots, (C_{k}, I_{k}))$ with $C_{1}, \ldots, C_{k}$ connected top-cycle-free diagrams and $\{I_{1}, \ldots, I_{k}\}$ a partition of $[1, j]$ into intervals with $1 \leqslant |I_{\ell}| \leqslant t_{1}(C_{\ell})$ for all $\ell$ and $j = |I_{1}| + \cdots + |I_{k}|$, and write $C = \beta((C_{1}, I_{1}), \ldots, (C_{k}, I_{k}))$ and let $c$ be the first terminal chord of $C$. Suppose $C$ contains a top cycle subdiagram $T$. Clearly $C - c$ is top-cycle-free since its connected components are exactly $C_{1}, \ldots, C_{k}$. Then it follows that $T$ must contain $c$ and, since cycles are 2-connected, be a subdiagram of $C_{\ell} \cup \{c\}$ for some $\ell$. Write $T = \{t_{1}, \ldots, t_{m}, c\}$, where $t_{1} < t_{2} < \cdots < t_{m} < c$ in the standard order. Without loss of generality we may assume that $T$ is chosen with $t_{2}$ maximum (in the standard order). Then, since $C_{\ell}$ is top-cycle-free, there is no path nested under $t_{2}$ from a chord of $T' = T - \{t_{1}, t_{2}, c\}$ to a right neighbor of $t_{2}$, implying that $T'$ is a subset of the source-sink group of $t_{2}$. By construction it follows that $t_{m}$ cannot cross $c$, a contradiction. So $C$ is top-cycle-free.

\begin{figure}[t]
\begin{center}
\includegraphics[scale=0.8]{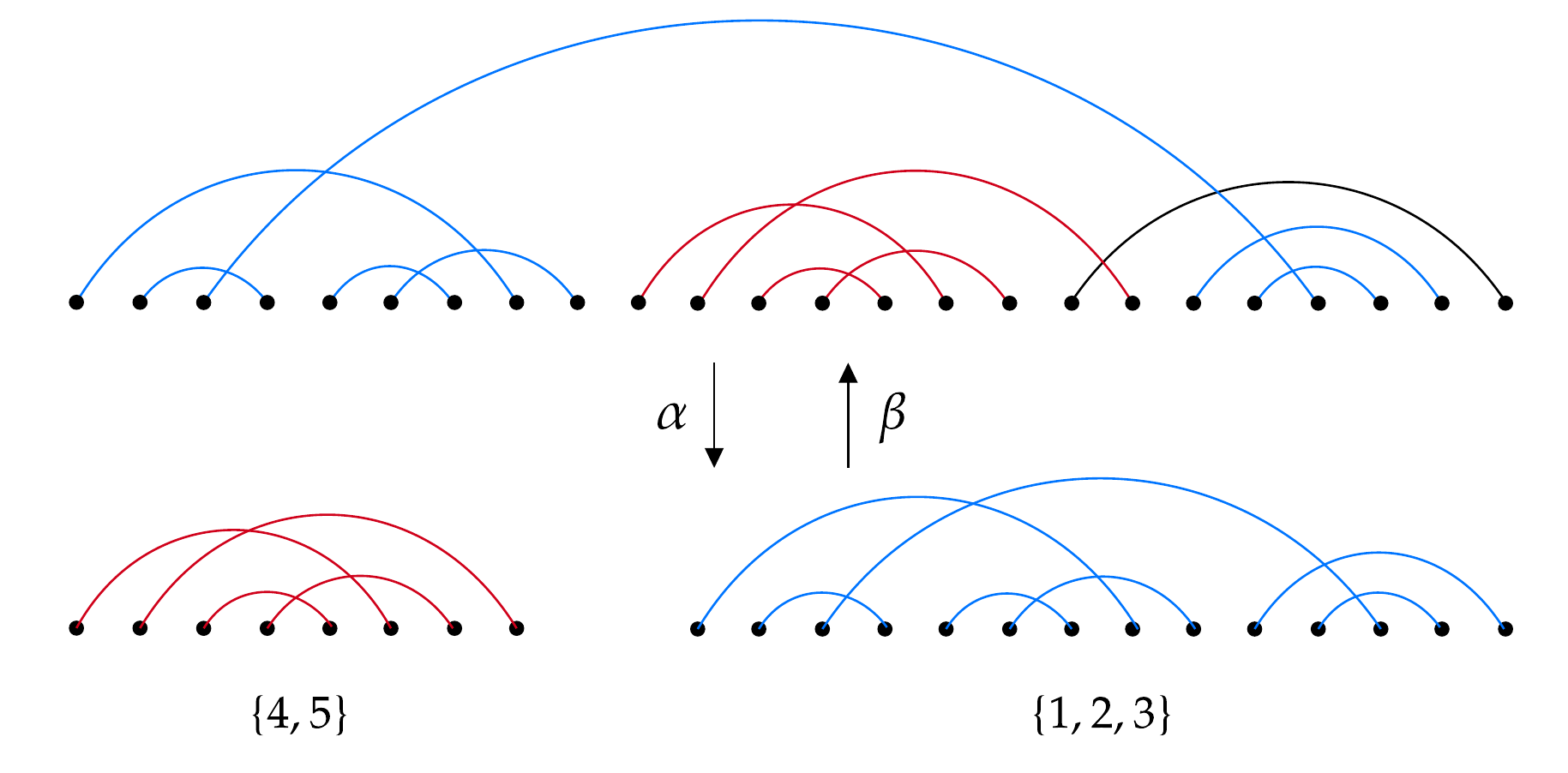}
\caption{A top-cycle-free diagram and its decomposition defining the maps $\alpha$ and $\beta$. The source-sink groups used in the construction of $\beta$ are indicated by the brackets below the two diagrams on the right side of the figure.}
\label{alpha_beta_decomp_fig2}
\end{center}
\end{figure}

Since $\alpha$ is a bijection and there is exactly one partition of $[1, j]$ into a fixed number of intervals with specified ordered sizes, this completes the proof (\ref{pow_identity_eq}) and Theorem~\ref{equation_sol_thm}. 



\section{A differential equation for the binomial 1-cocycle property}
\label{diff_eq_sec}

After seeing these results, a natural question immediately arises: why are connected chord diagrams the objects indexed
by these generating functions? In this section we present one attempt at an answer to this question, the roots of which go back to the work of Stein from the 1970s. The number $c_{n}$ of connected chord diagrams of size $n$ satisfies the well known recurrence 
\begin{align*}
c_{n+1} = \sum_{k=1}^{n}(2k - 1)c_{k}c_{n+1-k}
\end{align*}
due to Stein \cite{Stein1978a}. This corresponds to the differential equation 
\begin{align*}
C(x) - x = C(x)\left(2x{d \over dx} - 1\right)C(x)
\end{align*}
for the ordinary generating function $C(x)$ of connected diagrams. A bivariate version of this differential equation, known as the {\em renormalization group equation} in a physics context, turns out to be equivalent to the binomial 1-cocycle property when $\phi(z) = 1/(1-z)$, that is, all the weights are 1. Write $g_{i}(x) = i![y^{i}]G(x, y)$. 

\begin{theorem}
\label{diff_eq_thm}
If $G$ solves the functional equation 
\begin{align}
\label{gen_func_eq3}
G(x, y) = xL(\phi(G(x, y)))
\end{align}
with $\phi(z) = 1/(1 - z)$, then 
\begin{align*}
{\partial \over \partial y}G(x, y) = g_{1}(x)\left(2x{\partial \over \partial x} - 1\right)G(x, y) \quad \iff \quad L = L_{bin}.
\end{align*}
\end{theorem}

Before proving this, we require a lemma which essentially states that if $G$ is the weighted generating function over connected chord diagrams given by Theorem~\ref{equation_sol_thm} then it satisfies the differential equation. This lemma was originally proven by Marie and Yeats \cite{Marie2013} but we reprove it here for completeness. 

\begin{figure}[t]
\begin{center}
\includegraphics[scale=0.8]{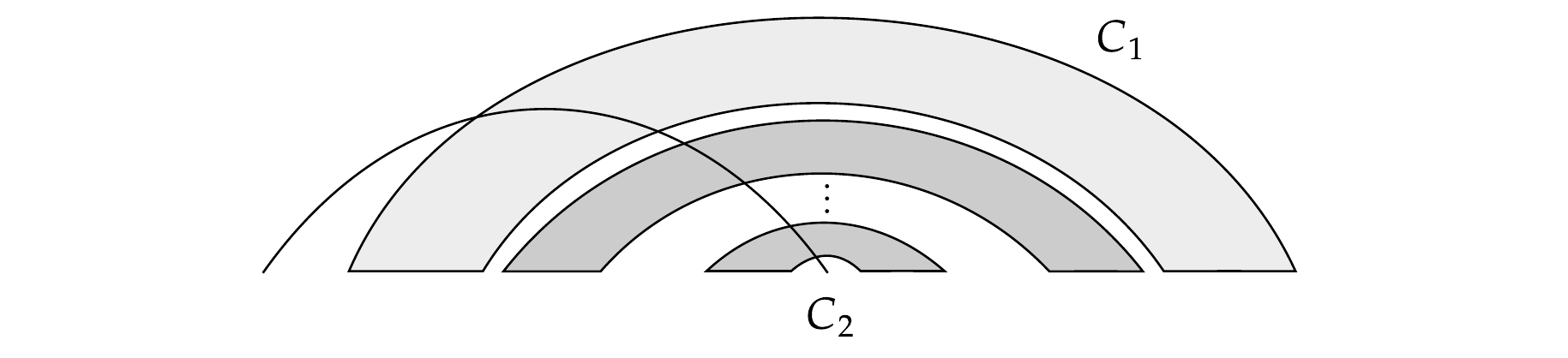}
\caption{The root-share decomposition $(C_{1}, C_{2})$ of $C$, with the insertion index suppressed.}
\label{root_share_fig}
\end{center}
\end{figure}

\begin{lemma}
\label{root_share_identity_lem}
For $n \geqslant 2$, we have
\begin{align*}
\sum_{\substack{C \in \mathcal{C} \\ |C| = n \\ t_{1}(C) \geqslant i}}f_{t_{1}(C)-i}f_{C} = \sum_{m = 1}^{n-i+1}(2(n - m) - 1)\sum_{\substack{C_{1} \in \mathcal{C} \\ |C_{1}| = m \\ t_{1}(C_{1}) \geqslant 1}}f_{t_{1}(C_{1})-1}f_{C_{1}}\sum_{\substack{C_{2} \in \mathcal{C} \\ |C_{2}| = n-m \\ t_{1}(C_{2}) \geqslant i-1}}f_{t_{1}(C_{2})-i+1}f_{C_{2}}.
\end{align*}
\end{lemma}
\begin{proof}
Let $C$ be a connected diagram of size $n$ with $t_{1}(C) \geqslant i$. Removing the root chord $c_{1}$ of $C$ gives an indecomposable diagram $C'$; let $C_{2}$ be the outermost component of $C'$, $C_{1} = C - C_{2}$, and $i$ be the endpoint of $C_{2}$ immediately left of the leftmost endpoint of $C_{1} - \{c_{1}\}$ in $C$. Then $(C_{1}, C_{2}, i)$ defines what has been called the {\em root-share decomposition} of $C$ (see Figure~\ref{root_share_fig}). Clearly we can recover $C$ from $(C_{1}, C_{2}, i)$. By construction, $C_{1}$ and $C_{2}$ are connected. Furthermore, by definition, $t_{1}(C_{2}) = t_{1}(C) - 1 \geqslant i - 1$ and, setting $m = |C_{1}|$, $i$ can take any value in $[1, 2(n-m)-1]$. Since each terminal chord of $C$ corresponds uniquely to a terminal chord of $C_{1} - \{c_{1}\}$ or $C_{2}$, it easily follows that 
\begin{align*}
f_{C} = f_{C_{1}}f_{t_{1}(C_{1})-1}f_{C_{2}}. 
\end{align*}
By all of the above observations, we infer the desired identity. 
\end{proof}

\begin{proof}[Proof of Theorem~\ref{diff_eq_thm}] 
We begin with the ``only if" direction. Suppose $G$ satisfies the tree-like equation (\ref{gen_func_eq3}) and 
\begin{align}
\label{rg_diff_eq}
{\partial \over \partial y}G(x, y) = g_{1}(x)\left(2x{\partial \over \partial x} - 1\right)G(x, y).
\end{align}
Write $f_{i} = {[y]L(y^{i}) \over i!}$ and $G(x, y) = \sum_{n \geqslant 1}\left(\sum_{i = 1}^{n}g_{i,n}{y^{i} \over i!}\right)x^{n}$. We proceed by induction on $n$ to simultaneously show that 
\begin{align*}
g_{i,n} = \sum_{\substack{C \in \mathcal{C} \\ |C| = n \\ t_{1}(C) \geqslant i}}f_{t_{1}(C)-i}f_{C} \quad \textrm{and} \quad L(y^{n}) = n!\sum_{i=1}^{n+1}f_{n+1-i}{y^{i} \over i!},
\end{align*}
that is, $L = L_{bin}$. Although proving that $G$ has the form given by Theorem~\ref{equation_sol_thm} is not part of the result statement, it seems necessary to generate a straightforward argument. Where before we began by expanding (\ref{gen_func_eq3}) only in $x$, we now expand in both $x$ and $y$ and then extract coefficients in $x$ to turn the functional equation into the identity
\begin{align}
\label{identity_yvar_eq}
\sum_{i=1}^{n+1}g_{i,n+1}{y^{i} \over i!} = \sum_{j=1}^{n}\Bigg(\sum_{k=1}^{n}\sum_{\substack{n_{1} + \cdots + n_{k} = n \\ n_{\ell} \geqslant 1}}\sum_{\substack{i_{1} + \cdots + i_{k} = j \\ 1 \leqslant i_{\ell} \leqslant m_{\ell}}}{j \choose i_{1}, \ldots, i_{k}}g_{i_{1},n_{1}} \cdots g_{i_{k},n_{k}}\Bigg){L(y^{j}) \over j!},
\end{align}
which holds for all $n \geqslant 1$. Applying the induction hypotheses and (\ref{bin_identity_eq}), we have
\begin{align*}
\sum_{i=1}^{n+1}g_{i,n+1}{y^{i} \over i!} &= \sum_{j=1}^{n}\Bigg(\sum_{k=1}^{n}\sum_{\substack{n_{1} + \cdots + n_{k} = n \\ n_{\ell} \geqslant 1}}\sum_{\substack{i_{1} + \cdots + i_{k} = j \\ 1 \leqslant i_{\ell} \leqslant m_{\ell}}}{j \choose i_{1}, \ldots, i_{k}} \\
&\qquad\qquad \times \sum_{\substack{C_{1} \in \mathcal{C} \\ |C_{1}| = n_{1} \\ t_{1}(C_{1}) \geqslant i_{1}}}f_{t_{1}(C_{1})-i_{1}}f_{C_{1}} \cdots \sum_{\substack{C_{k} \in \mathcal{C} \\ |C_{k}| = n_{k} \\ t_{1}(C_{k}) \geqslant i_{k}}}f_{t_{1}(C_{k})-i_{k}}f_{C_{k}}\Bigg)\sum_{i=1}^{j+1}f_{j+1-i}{y^{i} \over i!} \\
&= \sum_{j=1}^{n}\Bigg(\sum_{\substack{C \in \mathcal{C} \\ |C| = n+1 \\ t_{1}(C) = j+1}}f_{C}\Bigg)\sum_{i=1}^{j+1}f_{j+1-i}{y^{i} \over i!}.
\end{align*}
Then upon rearranging we obtain
\begin{align*}
\sum_{i=1}^{n+1}g_{i,n+1}{y^{i} \over i!} = \sum_{j=1}^{n+1}\Bigg(\sum_{\substack{C \in \mathcal{C} \\ |C| = n+1 \\ t_{1}(C) \geqslant i}}f_{t_{1}(C)-i}f_{C}\Bigg){y^{i} \over i!},
\end{align*}
so extracting coefficients gives the required expression for $g_{i,n+1}$. To derive $L(y^{n+1})$, we begin by noting that differential equation (\ref{rg_diff_eq}) expands into the identity
\begin{align}
\label{rge_identity_eq}
g_{i, n} = \sum_{m = 1}^{n-i+1}(2(n - m) - 1)g_{1,m}g_{i-1,n-m} \quad \textrm{for }  i \geqslant 2. 
\end{align}
Then 
\begin{align}
\nonumber \sum_{i=1}^{n+2}g_{i,n+2}{y^{i} \over i!} &= \sum_{i=2}^{n+2}\left(\sum_{m = 1}^{n-i+3}(2(n + 2 - m) - 1)g_{1,m}g_{i-1,n+2-m}\right){y^{i} \over i!} + g_{1,n+2}y \\
\nonumber &= \sum_{i=2}^{n+2}\Bigg(\sum_{m = 1}^{n-i+3}(2(n + 2 - m) - 1)\sum_{\substack{C_{1} \in \mathcal{C_{1}} \\ |C_{1}| = m \\ t_{1}(C_{1}) \geqslant 1}}f_{t_{1}(C_{1})-i}f_{C_{1}} \\
\nonumber &\qquad\qquad\qquad \times \sum_{\substack{C_{2} \in \mathcal{C} \\ |C_{2}| = n+2-m \\ t_{1}(C_{2}) \geqslant i-1}}f_{t_{1}(C_{2})-i+1}f_{C_{2}}\Bigg){y^{i} \over i!} + g_{1,n+2}y \\
\label{rhs_eq} &= \sum_{i=2}^{n+2}\Bigg(\sum_{\substack{C \in \mathcal{C} \\ |C| = n+2 \\ t_{1}(C) \geqslant i}}f_{t_{1}(C)-i}f_{C}\Bigg){y^{i} \over i!} + g_{1,n+2}y
\end{align}
by Lemma~\ref{root_share_identity_lem}. On the other hand, applying (\ref{identity_yvar_eq}) and (\ref{bin_identity_eq}) again we get
\begin{align}
\nonumber \sum_{i=1}^{n+2}g_{i,n+2}{y^{i} \over i!} &= \sum_{j=1}^{n+1}\Bigg(\sum_{k=1}^{n+1}\sum_{\substack{n_{1} + \cdots + n_{k} = n+1 \\ n_{\ell} \geqslant 1}}\sum_{\substack{i_{1} + \cdots + i_{k} = j \\ 1 \leqslant i_{\ell} \leqslant m_{\ell}}}{j \choose i_{1}, \ldots, i_{k}}g_{i_{1},n_{1}} \cdots g_{i_{k},n_{k}}\Bigg){L(y^{j}) \over j!} \\
\nonumber &= \sum_{j=1}^{n+1}\Bigg(\sum_{k=1}^{n+1}\sum_{\substack{n_{1} + \cdots + n_{k} = n+1 \\ n_{\ell} \geqslant 1}}\sum_{\substack{i_{1} + \cdots + i_{k} = j \\ 1 \leqslant i_{\ell} \leqslant m_{\ell}}}{j \choose i_{1}, \ldots, i_{k}} \\
\nonumber &\qquad\qquad \times \sum_{\substack{C_{1} \in \mathcal{C} \\ |C_{1}| = n_{1} \\ t_{1}(C_{1}) \geqslant i_{1}}}f_{t_{1}(C_{1})-i_{1}}f_{C_{1}} \cdots \sum_{\substack{C_{k} \in \mathcal{C} \\ |C_{k}| = n_{k} \\ t_{1}(C_{k}) \geqslant i_{k}}}f_{t_{1}(C_{k})-i_{k}}f_{C_{k}}\Bigg){L(y^{j}) \over j!} \\
\nonumber &= \sum_{j=1}^{n}\Bigg(\sum_{\substack{C \in \mathcal{C} \\ |C| = n+2 \\ t_{1}(C) = j+1}}f_{C}\Bigg)\sum_{i=1}^{j+1}f_{j+1-i}{y^{i} \over i!} + \sum_{\substack{C \in \mathcal{C} \\ |C| = n+2 \\ t_{1}(C) = n+2}}f_{C}{f(y^{n+1}) \over (n+1)!} \\
\label{lhs_eq} &= \sum_{j=1}^{n+1}\Bigg(\sum_{\substack{C \in \mathcal{C} \\ |C| = n+2 \\ i \leqslant t_{1}(C) \leqslant n+1}}f_{t_{1}(C)-i}f_{C}\Bigg){y^{i} \over i!} + \sum_{\substack{C \in \mathcal{C} \\ |C| = n+2 \\ t_{1}(C) = n+2}}f_{C}{f(y^{n+1}) \over (n+1)!}.
\end{align}
Comparing coefficients of (\ref{lhs_eq}) and (\ref{rhs_eq}), it follows that 
\begin{align*}
[y^{i}]f(y^{n+1}) = {f_{n+2-i} \over i!}
\end{align*}
for $i \geqslant 2$, as required. We conclude the proof by noting that the ``if" direction is directly obtained by combining (\ref{rge_identity_eq}), Lemma~\ref{root_share_identity_lem}, and Theorem~\ref{equation_sol_thm}. 
\end{proof}

While the ``if" direction of this theorem was proven in \cite{Marie2013}, the ``only if" direction is new. Since in general the weighted count of weighted chord diagrams is only double factorials when the weights are all 1, it is natural for this differential equation to only appear when $\phi(z) = 1/(1 - z)$. It is unclear whether analogous differential equations can be obtained in the case of arbitrary weights or as an equivalent condition to the divided power 1-cocycle property. 



\section{Diagram combinatorics from the generating functions}
\label{comb_sec}

We now turn to considering various bijective combinatorics arising from the generating functions obtained in Theorem~\ref{equation_sol_thm}.

\subsection{Combinatorial maps and triangulations of a disk}
\label{comb_map_triangle_sec}

We solved the tree-like equation
\begin{align}
\label{gen_func_eq4}
G(x, y) = xL(\phi(G(x, y)))
\end{align}
with weighted generating functions for certain classes of connected chord diagrams. What other combinatorial objects index these generating functions? Originally, in \cite{Marie2013}, Marie and Yeats solved the unweighted version of (\ref{gen_func_eq4}) with $L = L_{bin}$ by passing from chord diagrams to rooted plane binary trees. This essentially amounts to showing that there is a recursively-defined bijection, based on the root-share decomposition, between rooted connected chord diagrams and rooted plane binary trees with two recursively-defined properties. These properties are somewhat technical and based on a labeling of the leaves induced via the bijection by the intersection order on the chords of the associated diagram. 

Later, Courtiel, Yeats, and Zeilberger \cite{Courtiel2017b} discovered a much nicer type of object, combinatorial maps, that could substitute for connected chord diagrams. A {\em rooted combinatorial map} is a transitive permutation representation of the group $\langle \sigma, \alpha \mid \alpha^{2} = 1\rangle$ with a distinguished fixed point for the action of $\alpha$ (the root). It can be represented by a connected graph of half edges, each paired with at most one other half-edge to form an edge, with a cyclic order of the half-edges incident to each vertex and the root a vertex attached to a distinguished ``dangling" half-edge. Such a map inherits certain properties of this underlying graph; in particular, it is {\em bridgeless} if its graph is bridgeless, that is, 2-edge-connected. Other properties require a map-specific definition: it is {\em planar} if its Euler characteristic is 2 (see e.g. \cite{Courtiel2017b} for the definition of the Euler characteristic). 

\begin{theorem}[Courtiel, Yeats, and Zeilberger \cite{Courtiel2017b}]
\label{bridgeless_thm}
There exists a bijection $\theta$ between connected diagrams and bridgeless maps such that
\begin{itemize}
\item chords correspond to edges;
\item terminal chords correspond to vertices; 
\item the position $t_{1}$ of the first terminal chord corresponds to the indegree $\deg^{-}_{DFS, 1}$ of the root vertex $v_{1}$ under the orientation induced by the rightmost depth-first search. 
\end{itemize}
\end{theorem}

For brevity, we omit a definition of the rightmost depth-first search and associated orientation referenced in this theorem; see Section~5.5 of \cite{Courtiel2017b} for these details. This bijection $\theta$ also transfers the parameters involved in defining the weights $f_{C}$ and $\phi_{C}$, namely, the differences $t_{j} - t_{j-1}$ between the indices of the $j$th and $(j-1)$th terminal chords in the intersection order as well as chord valencies, into corresponding parameters on the associated bridgeless map. 

It turns out top-cycle-free diagrams also independently arise in the work of Courtiel, Yeats, and Zeilberger as the image of the the bijection $\theta$ on a natural subset of bridgeless maps: those which are planar. 

\begin{theorem}[Courtiel, Yeats, and Zeilberger \cite{Courtiel2017b}]
\label{planar_map_diagram_thm}
Under $\theta$, planar bridgeless maps are in bijection with connected top-cycle-free diagrams. 
\end{theorem}

Furthermore, changing ``connected" to ``indecomposable" corresponds to dropping ``bridgeless" for both this result and Theorem~\ref{bridgeless_thm}. It follows from the work of Courtiel, Yeats, and Zeilberger that just as there is a bridgeless map view of (\ref{binom_sol_eq}), there is also a planar bridgeless map view of (\ref{power_sol_eq}).  Several questions arise from these results, most notably: are there other well-behaved classes of combinatorial objects in bijection with top-cycle-free diagrams, especially classes that naturally index the generating function solutions to (\ref{gen_func_eq4})? Going further, can we explicitly enumerate such classes of objects? Although an affirmative answer to the latter question seems too much to hope for in the case of all connected diagrams and the binomial generating function, we can answer ``yes" to both questions for top-cycle-free diagrams and the divided power generating function.  

In this context, a {\em triangulation} $T$ is a plane graph in which every bounded face is a triangle; for technical reasons we also require that a single edge counts as a triangulation. Triangulation $T$ is {\em rooted} if it has a distinguished edge, the {\em root}; we will only work with rooted triangulations. As is standard in the literature, we consider such triangulations up to root-preserving isomorphism. {\em Exterior vertices} of $T$ are incident to its boundary face, while all other vertices of $T$ are {\em interior vertices}. 

\begin{theorem}
\label{diagram_triangle_thm}
There is a bijection $\omega$ between connected top-cycle-free diagrams with $n$ chords and $t_{1} = i$ and triangulations with $n - i$ interior vertices and $i+1$ exterior vertices. 
\end{theorem}

In 1964, William G. Brown \cite{Brown1964} explicitly enumerated (rooted) triangulations with $n$ interior vertices and $m + 3$ exterior vertices by deriving and solving a functional equation for the associated bivariate ordinary generating function, showing their number to be
\begin{align*}
{2(2m + 3)!(4n + 2m + 1)! \over m!(m + 2)!n!(3n + 2m + 3)!}.
\end{align*}
From Theorems \ref{planar_map_diagram_thm} and \ref{diagram_triangle_thm} it follows that we get an explicit count for the corresponding tree diagrams and planar bridgeless maps.

\begin{corollary}
\label{diagram_map_count_cor}
The number of connected top-cycle-free diagrams with $n$ chords and $t_{1} = i$ and planar bridgeless maps with $n$ edges and $\deg^{-}_{DFS, 1} = i$ is
\begin{align*}
{2(2i-1)!(4n-2i-3)! \over (i-2)!i!(n-i)!(3n-i-1)!} = {i-1 \over (4n - 2i - 1)(2n - i - 1)}{2i - 1 \choose i}{4n - 2i - 1 \choose n - i}.
\end{align*} 
\end{corollary}

\begin{figure}[t]
\begin{center}
\includegraphics[scale=0.8]{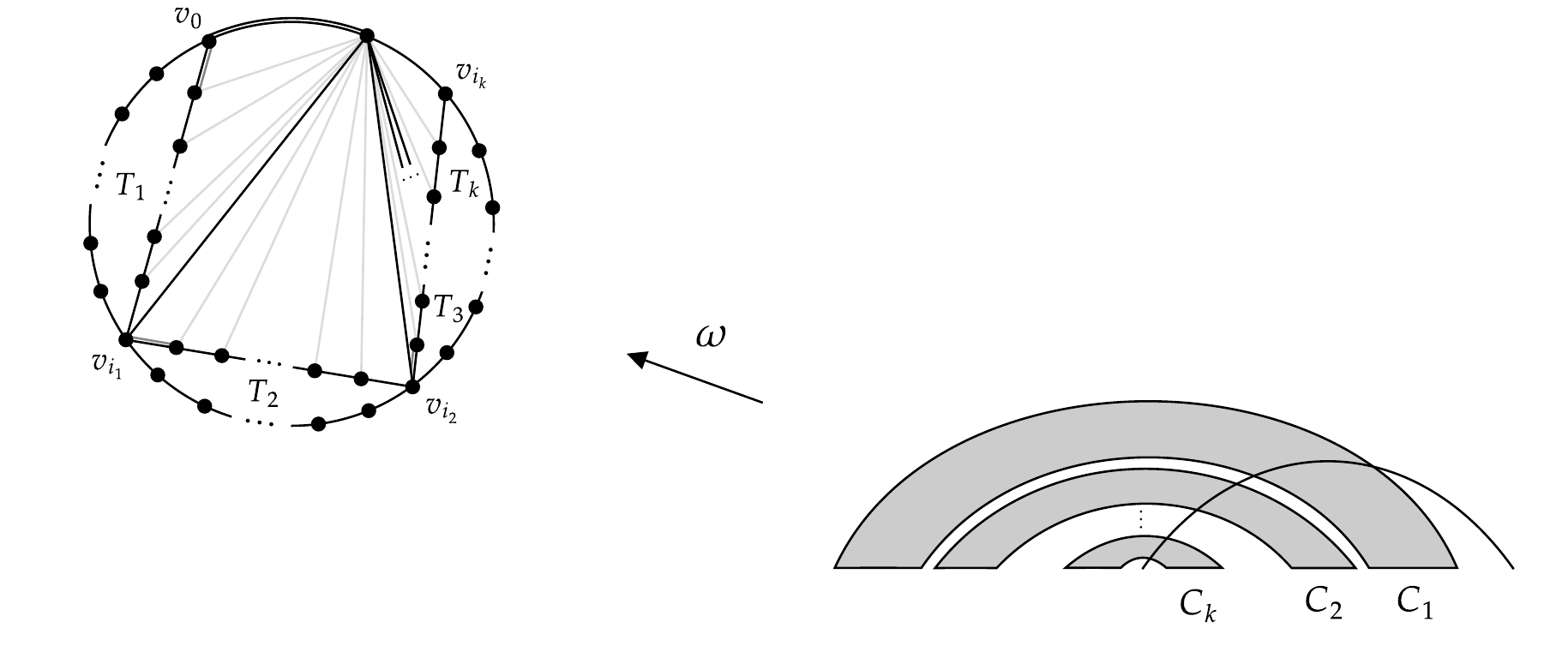}
\caption{The recursive construction of a triangulation and its image chord diagram under the bijection $\omega$. The exterior vertex labels are taken from the decomposition map $\gamma$.}
\label{triangulation_bij_fig}
\end{center}
\end{figure}

\begin{proof}[Proof of Theorem~\ref{diagram_triangle_thm}] 
We will use the unique decomposition described in sections \ref{bin_eq_sec} and \ref{pow_eq_sec} to recursively define our bijection. There is exactly one (connected top-cycle-free) diagram of size 1, namely, a single chord, and we map it to the unique (rooted) triangulation with no interior vertices and 2 exterior vertices: a single edge. So, explicitly, $\omega$ maps a single chord to a single edge. Now fix a connected top-cycle-free diagram $C$ of size $n \geqslant 2$. Applying the unique decomposition from previous sections, we get connected top-cycle-free subdiagrams $C_{1}, \ldots, C_{k}$ and integers $i_{1}, \ldots, i_{k}$ with $i_{1} + \cdots + i_{k} = t_{1}(C)-1$ and $1 \leqslant i_{\ell} \leqslant t_{1}(C_{\ell})$ for all $1 \leqslant \ell \leqslant k$; we have skipped the intervals $I_{1}, \ldots, I_{k}$ because their lengths $i_{1}, \ldots, i_{k}$ are sufficient to specify the decomposition and define the bijection $\omega$. Then inductively we obtain image triangulations $T_{1}, \ldots, T_{k}$ of $C_{1}, \ldots, C_{k}$ under $\omega$; write $t_{\ell} = t_{1}(C_{\ell})$ and $v_{\ell,0}, \ldots, v_{\ell,t_{\ell}}$ for the exterior vertices of $T_{\ell}$ read counterclockwise starting at the leftmost vertex $v_{\ell,0}$ of the root edge (so, in particular, $v_{\ell,0}v_{\ell,t_{\ell}}$ is the root edge). Then we construct triangulation $T$ as follows:
\begin{enumerate}
\item Join $T_{1}, \ldots, T_{k}$ in that order by identifying $v_{\ell-1, i_{\ell}}$ and $v_{\ell, 0}$. 
\item Add a new vertex $v$ and connect it with an edge to $v_{1,0}$ and $v_{k,i_{k}}$; this creates a single new bounded face. 
\item Root the graph at edge $v_{1,0}v$. 
\item While keeping the graph simple, add an edge from $v$ to every vertex incident to the newly created bounded face. 
\end{enumerate}
See Figure~\ref{triangulation_bij_fig} for a visual representation of the construction. We then set $\omega(C) = T$. It is easy to see that $T$ is indeed a triangulation. Furthermore, the number of vertices of $T$ is 
\begin{align*}
|T_{1}| + \sum_{\ell=2}^{k}(|T_{\ell}| - 1) = |C_{1}| + 1 + \sum_{\ell=2}^{k}|C_{\ell}| = |C|
\end{align*}
while the number of external vertices of $T$ is
\begin{align*}
1 + \sum_{\ell=1}^{k}i_{\ell} + 1 = t_{1}(C) + 1. 
\end{align*}
We have thus defined $\omega$; in particular, if we write $\delta$ for the function sending $((T_{1}, i_{1}), \ldots, (T_{k}, i_{k}))$ to $T$ then $\omega = \delta \circ \omega \circ \alpha$ (where $\omega$ acts on a tuple of diagram-integer pairs in the obvious way). It remains only to show that $\omega$ is a bijection. To that purpose, it suffices to define a unique decomposition $\gamma$ of a triangulation $T$ reversing the above construction, that is, inverting $\delta$. This follows from the observation that
\begin{align*}
\beta \circ \omega^{-1} \circ \gamma = \alpha^{-1} \circ \omega^{-1} \circ \delta^{-1} = (\delta \circ \omega \circ \alpha)^{-1} = \omega^{-1},
\end{align*}  
so we automatically get a recursively-defined inverse of $\omega$. So, let $T$ be a triangulation and $v_{0}, v_{1}, \ldots, v_{i}$ be the exterior vertices of $T$ read counterclockwise with $v_{0}v_{i}$ the root edge. Again reading counterclockwise, vertex $v_{i}$ has two or more external neighbors $v_{0}, v_{i_{1}}, v_{i_{2}}, \ldots, v_{i_{k-1}}, v_{i_{k}}$. Then deleting all edges incident to $v_{i}$ leaves a sequence of unrooted triangulations $T_{1}, \ldots, T_{k}$ pairwise joined at vertices $v_{i_{1}}, \ldots, v_{i_{k-1}}$ (see Figure~\ref{triangulation_bij_fig}). Root each $T_{\ell}$ at the unique edge incident to $v_{i_{\ell-1}}$ which was previously incident to a bounded face, where $i_{0} = 0$. Then it readily follows that setting $\gamma(T) = ((T_{1}, i_{1}), \ldots, (T_{k}, i_{k}))$ gives the desired decomposition; $T_{\ell}$ has exactly $i_{\ell}+1$ exterior vertices and, furthermore, clearly $\delta(\gamma(T)) = T$. This concludes the proof.
\end{proof}

It was first proved by Tutte \cite{Tutte1962} that there are 
\begin{align}
\label{triangle_count}
{2 \over n(n+1)}{4n+1 \choose n-1}
\end{align}
triangulations with $n$ internal vertices and 3 external vertices, while Wash and Lehman \cite{Walsh1975} proved that this also counts the number of rooted bridgeless planar maps with $n$ edges. Bijections connecting these two sets of objects were later given by Wormald \cite{Wormald1980}, Fusy \cite{Fusy2010}, and Fang \cite{Fang2018}, with the former two obtained recursively and the latter directly. Both our bijection between connected top-cycle-free diagrams and triangulations and the bijection of Courtiel, Yeats, and Zeilberger between planar bridgeless maps and connected top-cycle-free diagrams are recursively-defined. Obtaining a direct bijection mapping connected top-cycle-free diagrams to these other objects remains open.

\subsection{Catalan intervals and uniquely sorted permutations}
\label{interval_sorted_sec}

\newcommand{\Int}{\textnormal{Int}}
\newcommand{\NC}{\textnormal{NC}}
\newcommand{\PC}{\textnormal{PC}}

There are more than a few other well studied objects in bijection with triangulations and, thereby, connected top-cycle-free diagrams. We now review these objects and associated work to motivate the next section introducing new enumerative links involving chord diagram classes defined by forbidding certain subdiagrams. 

A plane binary tree $B$ can be represented as either a leaf $\circ$ or an ordered pair of plane binary trees $(B_{1}, B_{2})$. Tamari \cite{Tamari1962} defined a partial order on the set $\mathcal{B}_{n}$ of plane binary trees with $n$ non-leaf vertices and proved that it specified a lattice, the $n$\textsuperscript{th} {\em Tamari lattice} $\mathcal{L}^{T}_{n}$, whose covering relation is defined as follows: a plane binary tree B containing a subtree $A = ((B_{1}, B_{2}), B_{3})$ is covered by the plane binary tree $B'$ obtained by replacing $A$ by $(B_{1}, (B_{2}, B_{3}))$; this operation is known as right rotation. We define an {\em interval} in a poset $P$ as a pair $(x, y)$ of elements of $P$ such that $x \leqslant y$ and let $\Int(P)$ denote the set of all intervals of $P$. Chapoton \cite{Chapoton2006} proved that the cardinality of $\Int(\mathcal{L}^{T}_{n})$ is (\ref{triangle_count}); later, Bernardi and Bonichon \cite{Bernardi2009} provided a bijection between Tamari intervals and triangulations explaining their common count. Along with the bijection between triangulations and bridgeless planar maps, Fang \cite{Fang2018} obtained a bijection between bridgeless planar maps and Tamari intervals; all of these bijections passed through so-called ``sticky trees". 

In 2020, Colin Defant \cite{Defant2020b} proved that Tamari intervals are in bijection with certain pattern-avoiding permutations. Define a map $s$ on the set $\mathcal{S}_{n}$ of permutations of $[n]$ recursively as follows: $s$ sends the empty permutation to itself and, using one-line notation, sends a permutation $\sigma \in S_{n}$ to $s(\sigma_{L})s(\sigma_{R})n$, where $\sigma = \sigma_{L}n\sigma_{R}$. This is West's stack-sorting map, a variant of the stack-sorting algorithm introduced by Knuth \cite{Knuth1973} and studied extensively in West's Ph.D. thesis \cite{West1990}. There has been considerable interest in this map, especially with regards to permutations $\sigma$ with positive {\em fertility}, that is, $|s^{-1}(\sigma)|$, the number of preimages of $\sigma$ under $s$. Such permutations are referred to as {\em sorted}, and they are {\em uniquely sorted} if the fertility is 1. Let $\mathcal{U}_{n}$ denote the set of uniquely sorted permutation in $\mathcal{S}_{n}$; it was proved in \cite{Defant2020a} that $U_{n}$ is empty if $n$ is even. Defant proved that $\Int(\mathcal{L}^{T}_{n})$ is in bijection with permutations of $\mathcal{U}_{2n+1}$ that avoid the patterns $132$ or $231$. This is analogous to the classic result of Knuth that $132$- and $231$-avoiding permutations of $\mathcal{S}_{n}$ are counted by the $n$\textsuperscript{th} Catalan number $C_{n} = {1 \over n + 1}{2n \choose n}$. The analogy breaks down in an interesting way when considering other patterns of length 3. Whereas the Catalan numbers also count arbitrary permutations avoiding any other length 3 patterns, counts associated to non-Tamari Catalan intervals arise when excluding such patterns from uniquely sorted permutations of length $2n+1$. 


A {\em Dyck path} of size $n$ is a lattice path of North and East steps from $(0, 0)$ to $(n, n)$. There is a natural partial order $\leqslant_{S}$ on the Dyck paths $\mathcal{D}_{n}$ of size $n$ defined by setting $\pi \leqslant_{S} \pi'$ if $\pi$ lies weakly below $\pi'$. This defines a distributive lattice $\mathcal{L}^{S}_{n}$ known as the $n$\textsuperscript{th} {\em Stanley lattice}. De Sainte-Catherine and Viennot \cite{Sainte-Catherine1986} proved that 
\begin{align}
\label{stanley_count}
C_{n}C_{n+2} - C_{n+1}^{2} = {6 \over (n+1)(n+1)^{2}(n+3)}{2n \choose n}{2n+2 \choose n+1}
\end{align}
is the number of intervals of $\mathcal{L}^{S}_{n}$. We can similarly define a natural partial order on the set of noncrossing partitions $\NC_{n}$ of the set $[n]$, which generalize noncrossing diagrams, by ordering by refinement.  Kreweras \cite{Kreweras1972} proved that this poset is also a lattice, the $n$\textsuperscript{th} {\em Kreweras lattice} $\mathcal{L}^{K}_{n}$, and that the number of intervals in this lattice is 
\begin{align}
\label{kreweras_count}
{1 \over 2n + 1}{3n \choose n} 
\end{align} 
Note that, among other things, this is also the number of ternary trees with $n$ non-leaf vertices. In \cite{Pallo2003}, Pallo introduced a poset on $\mathcal{B}_{n}$ refining the Tamari lattice where rotation is restricted to occurring on subtrees of the path from the root to the rightmost leaf. This is known as the $n$\textsuperscript{th} {\em Pallo comb poset} $\PC_{n}$, and its intervals are counted by the coefficients of the generating series
\begin{align}
\label{pallo_count}
C(xC(x)),
\end{align}
where $C(x) = {1 - \sqrt{1 - 4x} \over 2x}$ is the generating function of the Catalan numbers. Finally, there is the antichain $\mathcal{A}_{n}$ of Catalan objects whose intervals are trivially counted by the Catalan numbers. Since plane binary trees, Dyck paths, and noncrossing partitions are all counted by the Catalan numbers, these five posets are collectively known as {\em Catalan posets}. By using bijections between Catalan objects to define the three lattices on a common groundset (e.g. Dyck paths \cite{Bernardi2009}), it can be shown that the Stanley lattice $\mathcal{L}^{S}_{n}$ is an extension of the Tamari lattice $\mathcal{L}^{T}_{n}$ which itself is an extension of the Kreweras lattice $\mathcal{L}^{K}_{n}$. Each of these three Catalan lattices has been the subject of extensive study throughout mathematics, with interest and applications in enumerative combinatorics and a diversity of other disciplines. 

Defant \cite{Defant2020b} completely characterized the enumerative relationship between Catalan intervals and length 3 pattern-avoiding uniquely sorted permutations. Letting $\mathcal{U}_{n}(\tau^{(1)}, \ldots, \tau^{(r)})$ denote the set of uniquely sorted permutations in $\mathcal{U}_{n}$ that avoid the patterns $\tau^{(1)}, \ldots, \tau^{(r)}$, Defant found bijections between $U_{2n+1}(321)$ and $\Int(\mathcal{L}^{S}_{n})$, between $U_{2n+1}(312, 1342)$ and $\Int(\mathcal{L}^{K}_{n})$, between $U_{2n+1}(231, 4132)$ and $\Int(\PC_{n})$, and between each of $U_{2n+1}(321)$, $U_{2n+1}(132, 312)$, $U_{2n+1}(132, 312)$, $U_{2n+1}(231, 312)$, and $\Int(\mathcal{A}_{n})$, the intervals of the Catalan antichain $\mathcal{A}_{n}$.

\subsection{Forbidden graphical subdiagrams}
\label{forbid_subdiagram_sec}

For a fixed set $\mathcal{X}$ of diagrams, let $\mathcal{D}_{n}(\mathcal{X})$ be the set of chord diagrams of size $n$ with no subdiagram lying in $\mathcal{X}$ and set $\mathcal{C}_{n}(\mathcal{X}) = \mathcal{D}_{n}(\mathcal{X}) \cap \mathcal{C}$, the connected {\em $\mathcal{X}$-free} diagrams. Writing $T_{i}$ for the top cycle diagram of size $i$, note that $\mathcal{C}_{n}(\{\})$ is of course $\mathcal{C}_{n}$, the set of size $n$ connected diagrams, and $\mathcal{C}_{n}(T_{\geqslant 3})$ is exactly $\mathcal{C}_{top}$, the set of size $n$ connected top-cycle-free diagrams. 

\begin{figure}[t]
\begin{center}
\includegraphics[scale=0.8]{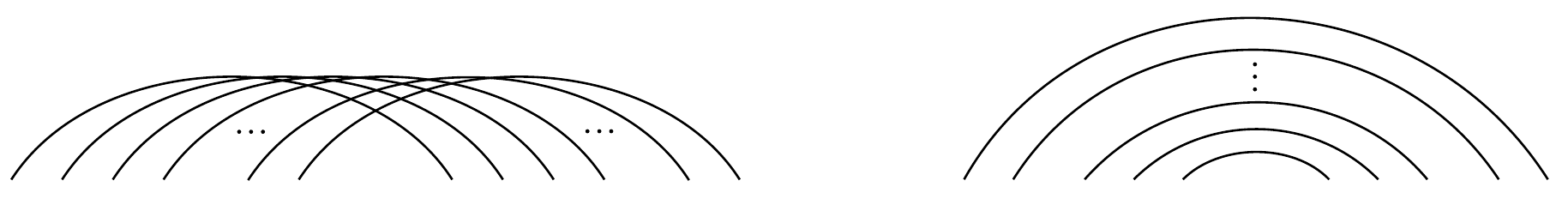}
\caption{The complete diagram and nesting diagram; the former is the unique representative of a complete graph.}
\label{complete_nesting_fig}
\end{center}
\end{figure}

Forbidding certain subdiagrams has previously been studied in the literature. Much past work has focused on forbidding complete subdiagrams $K_{k}$ and {\em nesting subdiagrams} $N_{k}$ of size $k$, traditionally referred to as {\em $k$-crossings} and {\em $k$-nestings} (see Figure~\ref{complete_nesting_fig}). It is well known that noncrossing diagrams $\mathcal{D}_{n}(K_{2})$ and nonnesting diagrams $\mathcal{D}_{n}(N_{2})$ are both counted by the Catalan numbers. Following work by Touchard \cite{Touchard1952}, Riordan \cite{Riordan1975}, and de Sainte-Catherine \cite{Sainte-Catherine1983} on the distribution of 2-crossings and 2-nestings, Gouyou-Beauchamps \cite{Gouyou-Beauchamps1989} studied the enumeration of $N_{3}$-free diagrams via involutions with no decreasing sequence of length 6, essentially giving a bijection between such diagrams and pairs of noncrossing Dyck paths. Chen et al. \cite{Chen2007} then extended these results to $k$-crossings and $k$-nestings, proving that their maximum size in a diagram defines a pair of symmetrically distributed statistics and therefore $|\mathcal{D}_{n}(K_{k})| = |\mathcal{D}_{n}(N_{k})|$. In particular, via a bijection between partitions and vacillating tableaux underlying their results they obtained the following. 

\begin{theorem}[Chen et al. \cite{Chen2007}]
The numbers of $K_{k}$-free and $N_{k}$-free diagrams of size $n$ are equal to the number of closed lattice paths of length $2n$ in the set 
\begin{align*}
\{(a_{1}, a_{2}, \ldots, a_{k-1}) \mid a_{1} \geqslant a_{2} \geqslant \cdots \geqslant a_{k-1} \geqslant 0, a_{i} \in \mathbb{Z}\}
\end{align*}
from the origin to itself with units steps in any coordinate direction or its negative. 
\end{theorem}

\begin{corollary}[Gouyou-Beauchamps \cite{Gouyou-Beauchamps1989}, Chen et al. \cite{Chen2007}]
\label{triangle_free_stanley_cor}
There is a bijection between $\mathcal{D}_{n}(K_{3})$ and $\Int(\mathcal{L}_{n}^{S})$, so both have cardinality $C_{n}C_{n+2} - C_{n+1}^{2}$.  
\end{corollary}

While the language we use here is inspired by graph theory and, in particular, the graph-theoretic nature of top cycle diagrams, classically in enumerative combinatorics the forbidding of substructures takes the form of pattern avoidance. Inspired by the broad literature on pattern avoidance in permutations (see \cite{Kitaev2011}), several authors have defined and studied pattern avoidance in set partitions and matchings, that is, chord diagrams \cite{Klazar1996, Jelinek2005, Sagan2010, Bloom2013}. Our definition of forbidding subdiagrams in particular matches the definitions of Jel\'{i}nek \cite{Jelinek2005} and Sagan \cite{Sagan2010} of avoiding matching patterns. Following Chen et al., Jel\'{i}nek considered forbidding other subdiagrams of size three: $D_{231} = \{(1, 5), (2, 6), (3, 4)\}$, $D_{312} = \{(1, 6), (2, 4), (3, 5)\}$, $D_{213} = \{(1, 5), (2, 4), (3, 6)\}$, and $D_{132} = \{(1, 4), (2, 6), (3, 5)\}$. These are {\em permutation diagrams}, so called because they are defined (and indexed) by the permutation determined by their sinks (see Section~\ref{psi_conn_sec} for a full definition). Jel\'{i}nek constructed bijections between $\mathcal{D}_{n}(K_{3})$ and $\mathcal{D}_{n}(D_{231})$ and between $\mathcal{D}_{n}(D_{213})$ and $\mathcal{D}_{n}(D_{132})$ preserving certain substructures. For the latter, he actually passed through top-cycle-free diagrams, giving bijections between both of those sets and $\mathcal{D}_{n}(T_{\geqslant 3})$. Bloom and Elizalde \cite{Bloom2013} later enumerated these classes of diagrams, simplifying prior bijections and deriving the algebraic generating function of $\mathcal{D}_{n}(D_{213})$ and, thereby, $\mathcal{D}_{n}(D_{132})$ and $\mathcal{D}_{n}(T_{\geqslant 3})$. 

The combinatorial significance of top cycle diagrams motivates considering other diagrams with graph-theoretic relevance and, in particular, studying sets of diagrams which forbid such graphical subdagrams. In graph theory, graphs forbidding certain cycles as induced subgraphs are of great interest. Most notably this includes triangle-free graphs (forbidding a cycle of size 3), forests (forbidding all cycles), chordal graphs (forbidding cycles of size $4$ or greater), and bipartite graphs (forbidding odd size cycles). By forbidding the corresponding top and bottom cycles, we get the chord diagram versions of each of these graph classes, namely $\mathcal{D}_{n}(T_{3})$, $\mathcal{D}_{n}(T_{\geqslant 3},B_{\geqslant 3})$, $\mathcal{D}_{n}(T_{\geqslant 4},B_{\geqslant 4})$, and $\mathcal{D}_{n}(\{T_{2k+1},B_{2k+1}\}_{k \geqslant 1})$, where $B_{i}$ is the bottom cycle diagram of size $i$.\footnote{Note that $B_{3} = T_{3}$.} In addition to the top-cycle-free diagrams $\mathcal{D}_{n}(T_{\geqslant 3})$, we naturally also consider bottom-cycle-free diagrams $\mathcal{D}_{n}(B_{\geqslant 3})$. As discussed previously, we already know that triangle-free diagrams $\mathcal{D}_{n}(T_{3})$ are in bijection with, among other things, intervals $\Int(\mathcal{L}^{S}_{n})$ of the Stanley lattice (Corollary~\ref{triangle_free_stanley_cor}). By Theorem~\ref{diagram_triangle_thm} and previously mentioned results of Bernardi and Bonichon \cite{Bernardi2009}, we also know that connected top-cycle-free diagrams $\mathcal{C}_{n}(T_{\geqslant 3})$ are equinumerous with another Catalan lattice. 

\begin{corollary}
\label{top_cycle_tamari_cor}
There is a bijection between $\mathcal{C}_{n}(T_{\geqslant 3})$ and $\Int(\mathcal{L}^{T}_{n})$, so both have cardinality ${2 \over n(n+1)}{4n+1 \choose n-1}$. 
\end{corollary}

The fact that the Stanley lattice extends the Tamari lattice then reflects the trivial observation that $C_{n}(T_{\geqslant 3})$ is a subset of $D_{n}(T_{3})$. If we further exclude bottom cycles, we get the tree diagrams $C_{n}(T_{\geqslant 3}, B_{\geqslant 3})$; comparing manual calculations of their count for small $n$ with OEIS sequence \href{https://oeis.org/A001764}{A001764} \cite{OEIS2021}, we conjecture the following. 

\begin{conjecture}
There is a bijection between $C_{n}(T_{\geqslant 3}, B_{\geqslant 3})$ and $\Int(\mathcal{L}^{K}_{n})$, so both have cardinality ${1 \over 2n+1}{3n \choose n}$. 
\end{conjecture}

It is notable that to pass from the Stanley case to the Tamari case we exclude an infinite set of cycle diagrams {\em and} assert connectedness, while to further pass to the Kreweras case we need only further exclude another infinite set of cycle diagrams. This may reflect a closer relationship between the Tamari and Kreweras lattices than between the Stanley and Tamari lattices, as is also suggested by e.g. the Dyck path realizations of these lattices given by Bernardi and Bonichon \cite{Bernardi2009}. 

Manual calculations similarly indicate that connected chordal diagrams $\mathcal{C}_{n}(T_{\geqslant 4}, B_{\geqslant 4})$, bipartite diagrams $\mathcal{C}_{n}(\{T_{2k+1}, B_{2k+1}\}_{k \geqslant 1})$, and bottom-cycle-free diagrams $\mathcal{C}_{n}(B_{\geqslant 3})$ are enumerated by simple, combinatorially significant sequences. 

\begin{conjecture}[\href{https://oeis.org/A001246}{A001246}]
The cardinality of $\mathcal{C}_{n}(T_{\geqslant 4}, B_{\geqslant 4})$ is $C_{n}^{2}$. 
\end{conjecture}

\begin{conjecture}[\href{https://oeis.org/A000264}{A000264}]
The cardinality of $\mathcal{C}_{n}(\{T_{2k+1}, B_{2k+1}\}_{k \geqslant 1})$ is the same as the set of all 3-edge-connected rooted cubic maps with $2n$ vertices and a distinguished Hamiltonian cycle. 
\end{conjecture}

\begin{conjecture}[\href{https://oeis.org/A064062}{A064062}]
\label{conn_bottom_cycle_conj}
The cardinality of $\mathcal{C}_{n}(B_{\geqslant 3})$ is 
\begin{align*}
C(2; n-1) = {1 \over n}\sum_{k=0}^{n-1}{2n \choose n - 1 - k}{n - 1 + k \choose k}.
\end{align*} 
\end{conjecture}

Neither of these three sets is equinumerous with the intervals of a known Catalan poset, although the possibility of such an association remains open for the latter two. The first conjecture clearly calls for a bijection between chordal diagrams of size $n$ and ordered pairs of Catalan objects of size $n$. With regards to the third conjecture, Bloom and Elizalde \cite{Bloom2013} also enumerated diagrams forbidding a pair of permutation diagrams of size $3$; in particular, they found that the {\em generalized Catalan numbers} $C(2; n)$ of (\ref{conn_bottom_cycle_conj}) also count $\{D_{213}, D_{132}\}$-free diagrams. This interestingly relates to the result of Jel\'{i}nek \cite{Jelinek2005} that top-cycle-free diagrams are in bijection with $D_{213}$-free diagrams and $D_{132}$-free diagrams. It also further motivates our additional conjecture, informed by numerical evidence, that $|\mathcal{C}_{n}(B_{\geqslant 3})| \leqslant |\mathcal{C}_{n}(T_{\geqslant 3})|$ for all $n$. 

All the enumerative conjectures we have made thus far have involved sets of connected diagrams. While none of their not-necessarily-connected versions have count sequences that appear on the OEIS, there is a classic functional equation relating the (ordinary) generating function $D(x)$ of chord diagrams and the generating function $C(x)$ of connected chord diagrams, 
\begin{align*}
D(x) = 1 + C(xD(x)^{2}). 
\end{align*}
Note that the empty diagram is not considered connected; see e.g. \cite{Flajolet2000}. This is obtained by decomposing a diagram by the component containing the root chord and viewing the rest of the diagram as attached to its chords. This straightforwardly generalizes to any class $\mathcal{D}(\mathcal{Z})$ determined by forbidding a fixed set $\mathcal{Z}$ of connected diagrams and its associated class of connected $\mathcal{Z}$-free diagrams $\mathcal{C}(\mathcal{Z})$. In particular, writing $G(x)$ and $F(x)$ for the generating functions of $\mathcal{D}(\mathcal{Z})$ and $\mathcal{C}(\mathcal{Z})$, we have 
\begin{align}
\label{conn_forbid_func_eq}
G(x) = 1 + F(xG(x)^{2}). 
\end{align}
With a routine expansion this translates into the recurrence 
\begin{align}
\label{conn_forbid_recurr_eq}
a_{n} = \sum_{k=1}^{n}b_{k}\sum_{n_{1} + \cdots n_{2k} = n - k}a_{n_{1}} \cdots a_{n_{2k}},
\end{align}
where $a_{n} = [x^{n}]G(x)$ and $b_{n} = [x^{n}]F(x)$. In certain cases it may be possible to apply a formula for $F(x)$ or $b_{n}$ and either (\ref{conn_forbid_func_eq}) or (\ref{conn_forbid_recurr_eq}) to obtain an explicit expression for $G(x)$ or $a_{n}$. Alternatively, it may feasible to use these equations to translate alternative combinatorial interpretations of $\mathcal{C}(\mathcal{Z})$ into a related combinatorial interpretation of $\mathcal{D}(\mathcal{Z})$. 

The conjectured appearance of connectedness in the enumerative study of forbidden subdiagrams suggests that other notions of connectivity, in particular indecomposability and 1-terminality, may also find significance in this context. Note that none of the three notions of connectivity has a forbidden subdiagram characterization. As we mentioned in Section~\ref{comb_map_triangle_sec}, indecomposable top-cycle-free diagrams are known to be in bijection with planar maps. We speculate that excluding other graphical diagrams may correspond to other natural types of combinatorial maps, but were not able to easily identify any such possible connections. On the other hand, applying 1-terminality seems to yield many enumerative links. Write $\mathcal{T}_{n}(\mathcal{X})$ for the set of $\mathcal{X}$-free 1-terminal diagrams of size $n$.

\begin{theorem}
The cardinality of $\mathcal{T}_{n}(T_{\geqslant 3}) = \mathcal{T}_{n}(T_{\geqslant 3}, B_{\geqslant 3})$ is $C_{n-1}$. 
\end{theorem} 

Note that the equality is nontrivial but straightforward. We prove this in Section~\ref{psi_map_sec} by exhibiting a bijection to noncrossing diagrams of size $n-1$. 

\begin{conjecture}[\href{https://oeis.org/A117106}{A117106}]
The cardinality of $\mathcal{T}_{n}(T_{3})$ is the same as the set of semi-Baxter permutations of length $n-1$,
\begin{align*}
{24 \over (n-2)(n-1)^{2}n(n+1)}\sum_{j=0}^{n-1}{n-1 \choose j+2}{n+1 \choose j}{n+j+1 \choose j+1}.
\end{align*} 	
\end{conjecture}

\begin{conjecture}[\href{https://oeis.org/A001181}{A001181}]
The cardinality of $\mathcal{T}_{n}(\{T_{2k+1},B_{2k+1}\}_{k \geqslant 1})$ is the same as the set of Baxter permutations of length $n-1$,
\begin{align}
\label{baxter_eq}
\sum_{k=1}^{n-1}{{n-2 \choose k-1}{n \choose k}{n \choose k+1} \over {n \choose 1}{n \choose 2}}.
\end{align}
\end{conjecture}

\begin{conjecture}[\href{https://oeis.org/A006318}{A006318}]
The cardinality of $\mathcal{T}_{n}(B_{\geqslant 3})$ is the $(n-2)$\textsuperscript{nd} Schr\"{o}der number, 
\begin{align*}
S_{n-2} = \sum_{k=0}^{n-2}C_{k}{n-2+k \choose n-2-k}. 
\end{align*}
\end{conjecture}

See the work of Chung et al. \cite{Chung1978} and Bouvel et al. \cite{Bouvel2018} for definitions of Baxter and semi-Baxter permutations and proofs of the closed-form expressions for their counts. Note that semi-Baxter permutations are a relaxation of Baxter permutations, matching the fact that bipartite diagrams are triangle-free. 

Clearly interesting enumerative structure seems to emerge from forbidding graphical subdiagrams and applying various notions of connectedness, warranting further study. In addition to cycles, excluding most notably complete graphs, trees, paths, and bipartite graphs has yielded structurally rich graph classes; forbidding their diagram representations may also find relevance in this context. Additionally, there are several natural questions that follow from Theorem~\ref{diagram_triangle_thm}. Are there bijections proving the above conjectures that translate between interesting statistics on the associated objects? For example, (\ref{baxter_eq}) is graded by the number of descents in the Baxter permutations, suggesting an associated statistic on 1-terminal bipartite diagrams. Finally, besides top-cycle-free diagrams does the inverse $\beta$ preserve any other excluded subdagram class of connected diagrams?












\section{The bijection $\psi$}
\label{psi_map_sec}

In order to understand the solution (\ref{binom_sol_eq}) to the binomial tree-like equation, Courtiel and Yeats \cite{Courtiel2017a} investigated the asymptotics of the parameters defining the solution, including the index of the first terminal chord, the number of terminal chords, and the differences between the indices of consecutive terminal chords in the intersection order. For the former they obtained a recurrence relation for the number $c_{n, k}$ of connected diagrams with $n$ chords such that $t_{1} \geqslant n - k$. In addition to recursively computing the associated exponential generating function and estimating the asymptotics of $c_{n, k}$, they used this result to derive the following count. 

\begin{theorem}[Courtiel and Yeats {\cite[Corollary 11]{Courtiel2017a}}]
\label{1term_count_thm}
The number of connected diagrams with $n$ chords and exactly one terminal chord is $(2n - 3)!!$. 
\end{theorem}

This is the same as the number of (not necessarily connected) diagrams with $n-1$ chords. Courtiel and Yeats used the root-share decomposition from Section \ref{diff_eq_sec} to inductively prove this count. In this section we describe and study a bijection $\psi: \mathcal{T}_{n} \rightarrow \mathcal{C}_{n-1}$ between the set $\mathcal{T}_{n}$ of 1-terminal diagrams with $n$ chords and the set $\mathcal{C}_{n-1}$ of diagrams with $n-1$ chords. This map was actually first discovered by Karen Yeats almost a decade ago in the course of the original work with Marie \cite{Marie2013} on the chord diagram solution to (\ref{bin_func_eq}), but never published or, to our knowledge, significantly studied.\footnote{Personal communication.} Yeats used a recursive formulation related to the inductive proof of Theorem~\ref{1term_count_thm}; here we describe and concentrate on a new, simpler formulation. We will return to Yeats' formulation in Section~\ref{factorial_objects_sec} in the context of Stirling permutations. 

\begin{figure}[t]
\begin{center}
\includegraphics[scale=0.74]{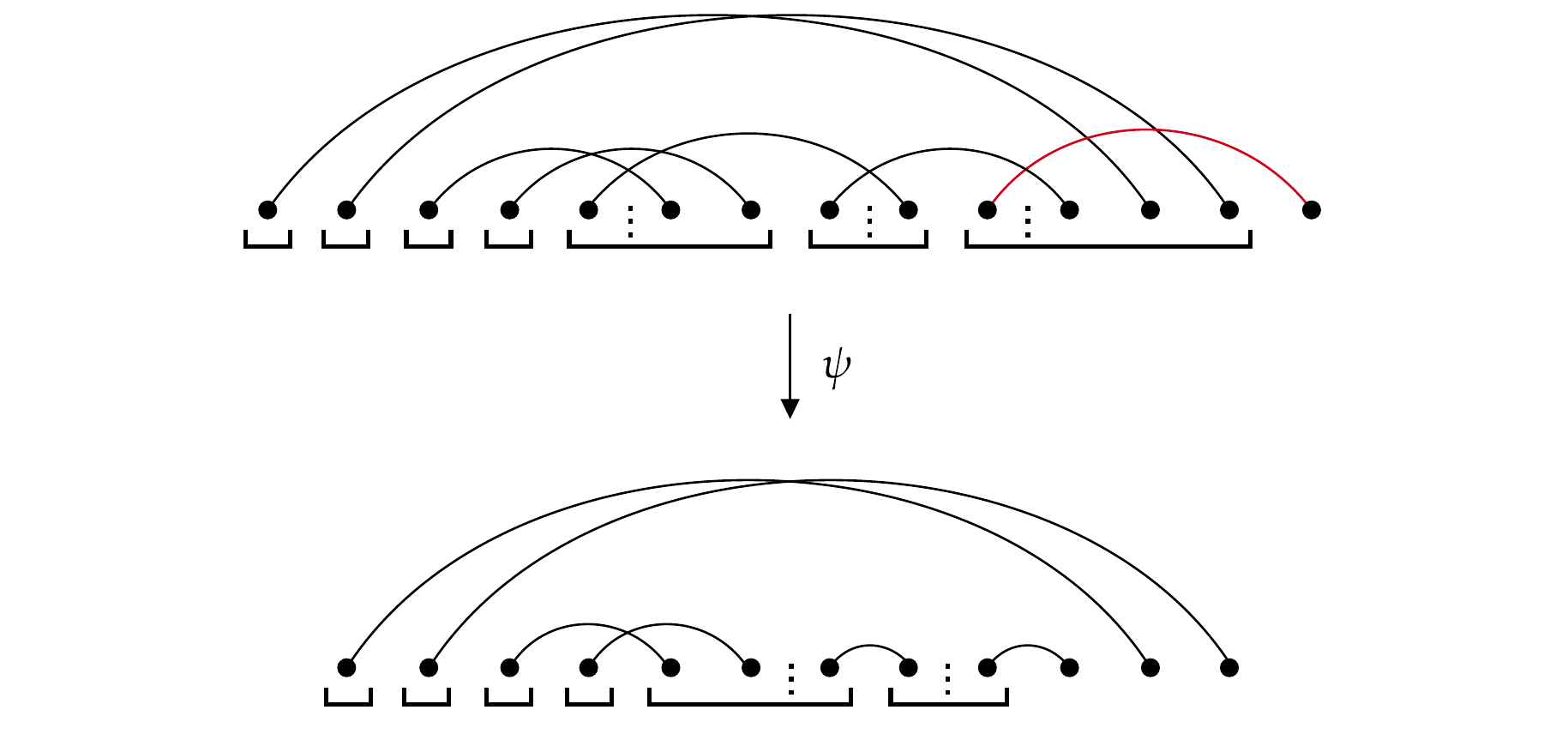}
\caption{An example of a 1-terminal diagram and its image diagram under the bijection $\psi$. The source-sink groups ``flipped" into sink-source groups by $\psi$ are indicated by the horizontal brackets, with the ``flip axes" indicated by the dotted lines, and the terminal chord of $T$ is indicated in red.}
\label{psi_map_fig}
\end{center}
\end{figure}

Consider $T \in \mathcal{T}_{n}$ with chords $c_{1} < c_{2} < \cdots < c_{n}$ in the standard order (which agrees with the intersection order by 1-terminality). Recall from Section~\ref{bin_eq_sec} that the endpoints of $T$ partition into a sequence of source-sink groups, one for each source of $T$. Since $T$ is 1-terminal, by Lemma~\ref{comp_nbhr_lem} these groups only contain the source and maximal interval of sinks immediately proceeding the source. Then we define $\psi(T)$ to be the diagram obtained by
\begin{enumerate}
\item[(1)] moving each source to the end of its source-sink group, and
\item[(2)] deleting the formerly terminal chord $c_{n}$. 
\end{enumerate}
Figure~\ref{psi_map_fig} displays a representative example of $T$ and $\psi(T)$. Note that the bijection $\psi$ induces a mapping, given by the construction, between the chords of $T$ and the chords of its image $\psi(T)$. We can therefore abuse notation and write $\psi(c)$ for the image of a non-terminal chord $c$ under this induced mapping. 

Step (1) converts each source-sink group in $T$ into a corresponding {\em sink-source group} in $\psi(T)$. Observe that we can also dually characterize this as uncrossing each non-terminal chord $c$ with its rightmost right neighbor $c'$ by moving the source of $c'$ immediately ahead of the sink of $c$ while maintaining the relative order of all other endpoints. In particular, the following key property holds. 

\begin{lemma}
\label{psi_right_nhbr_lem}
For $T \in \mathcal{T}_{n}$, non-terminal chord $c \in T$ has $k$ right neighbors if and only if $\psi(c)$ has $k-1$ right neighbors.  
\end{lemma}

The central claim is of course the following. 

\begin{theorem}
\label{psi_bijection_thm}
For all $n \geqslant 2$, the map $\psi$ is a bijection sending $\mathcal{T}_{n}$ to $\mathcal{C}_{n-1}$. 
\end{theorem}
\begin{proof}
Clearly the codomain of $\psi$ is $\mathcal{C}_{n-1}$, so it suffices to exhibit the inverse of $\psi$. To that end, let $C \in \mathcal{C}_{n-1}$ and define $\chi(C)$ to be the diagram obtained by
\begin{enumerate}
\item[(1')] concatenating a single chord to the end of $C$, and
\item[(2')] moving the source of each sink-source group to the beginning of its sink-source group. 
\end{enumerate}
Note that we could also skip using (2') on the terminal chord and instead of (1') directly add a chord covering the rightmost maximal interval of sinks. It is readily apparent that (1') and (2') invert (2) and (1), respectively, so $\chi$ is the inverse of $\psi$. 
\end{proof}

Reviewing Definition~\ref{traced_subdiagram_def} and the inverse map $\beta$ used to prove Theorem~\ref{equation_sol_thm}, we can see that through the use of source-sink groups $\psi$ is closely related to the decomposition of Section~\ref{bin_eq_sec}. In particular, we can apply $\psi$ to convert the traced-subdiagram-based decomposition of 1-terminal diagrams into a similar decomposition of arbitrary diagrams, where sink-source groups are used to define an alternate notion of traced subdiagrams. Although we leave the details out for brevity, this similarly gives a proof that the exponential generating function $C(x)$ of chord diagrams satisfies
\begin{align*}
C(x) = \int_{0}^{x} {dy \over 1 - C(y)}.
\end{align*} 

In the weighted generating functions (\ref{binom_sol_eq}) and (\ref{power_sol_eq}) 1-terminal diagrams index the terms on the diagonal where the exponent of $y$ equals the exponent of $x$. This raises a question: under $\psi$, to what do 1-terminal top-cycle-free diagrams map? The answer is illuminating.

\begin{lemma} 
\label{topcycle_1term_lem}
Let $T$ be top-cycle-free. Then $T$ is 1-terminal if and only if $T$ is a tree and every non-terminal chord has exactly one right neighbor. 
\end{lemma}
\begin{proof}
The ``if" direction holds by definition. For the other direction, assume $T$ is 1-terminal. Note that bottom cycle diagrams have chords with two right neighbors, so it suffices to show that every non-terminal chord has exactly one right neighbor (since $T$ consequently must be a tree). If not, then there exists non-terminal $c \in \mathcal{T}$ with at least two right neighbors $d$ and $d'$. Since $T$ is top-cycle-free $d$ and $d'$ do not cross and we may assume that $d$ is nested under $d'$. Choose $e$ as far right as possible such that there is a nonnesting $d$-$e$ path nested under $d'$. Then either $e$ has a right neighbor crossing $d'$, in which case $T$ contains a top cycle subdiagram, or $e$ is terminal. In either case we get a contradiction. 
\end{proof}

Note that we could also obtain this structural characterization of 1-terminal top-cycle-free diagrams using a quicker induction argument, but the above proof offers more insight. It in particular reflects the fact that bottom cycle diagrams are not 1-terminal. 

\begin{proposition} 
\label{psi_noncrossing_prop}
The map $\psi$ restricts to a bijection between 1-terminal top-cycle-free diagrams of size $n$ and noncrossing diagrams of size $n-1$. 
\end{proposition}
\begin{proof}
Let $T$ be top-cycle-free and 1-terminal. By Lemmas \ref{topcycle_1term_lem} and \ref{psi_right_nhbr_lem} each chord in the image $\psi(T)$ has no right neighbors, implying that there are no crossings at all. We can similarly infer that for each noncrossing diagram $C$ every non-terminal chord in $\psi^{-1}(C)$ has exactly one right neighbor, as required. 
\end{proof}

Among other things, this sheds further light on the consequence of Theorem~\ref{diagram_triangle_thm} that 1-terminal top-cycle-free diagrams are in bijection with triangulations with no interior vertices and, therefore, are counted by the Catalan numbers. In a nutshell, Lemma~\ref{topcycle_1term_lem} and Proposition~\ref{psi_noncrossing_prop} together imply that we can think of 1-terminal top-cycle-free diagrams as the connectivity 1 equivalent to noncrossing diagrams, since trees have connectivity 1. As indicated by Lemma~\ref{topcycle_1term_lem}, they are also minimally 1-terminal in the sense that each chord has the minimum number of right neighbors required to be 1-terminal. 

With these observations in mind we can generalize Proposition~\ref{psi_noncrossing_prop} to get a connectivity $k$ equivalent to noncrossing diagrams. Call a diagram $C$ {\em $k$-terminal-minimal} if all but its last $k-1$ chords have exactly $k$ right neighbors. Clearly $C$ is $k$-terminal and furthermore, in particular, 1-terminal top-cycle-free diagrams are 1-terminal-minimal. We will require a basic fact about connectivity, followed by two important statements about the relationship between $k$-connectivity, $k$-terminality, and $\psi$. For a simple graph $G$, write $G[A]$ for the induced subgraph on $A \subseteq V(G)$. 
 

\begin{lemma}
\label{graph_conn_lem}
If $G - A$ is $k$-connected, $G[A]$ is either $k$-connected or has size at most $k$, $\deg(V(G) - A) \geqslant \min\{k, |A|\}$, and $\deg(A) \geqslant k$, then $G$ is $k$-connected.  
\end{lemma}

We omit the straightforward proof of this fact. 

\begin{proposition}
\label{kterm_kconn_prop}
If diagram $C$ is $k$-terminal and has at least $k+1$ chords then it is $k$-connected. 
\end{proposition}
\begin{proof}
Let $c$ be the root chord of a $k$-terminal diagram $C$ of size at least $k+1$. By definition $C - c$ is $k$-terminal, so it is either the complete diagram of size $k$ or we may inductively assume that it is $k$-connected. In the former case $C$ is also complete and thus $k$-connected, while in the latter case the fact that $c$ has $k$ right neighbors in $C$ implies that $C$ is $k$-connected by Lemma~\ref{graph_conn_lem} with $A = \{c\}$. 
\end{proof}

\begin{proposition}
\label{psi_kterm_prop}
The map $\psi$ restricts to a bijection between $k$-terminal diagrams of size $n$ and $(k-1)$-terminal diagrams of size $n-1$. 
\end{proposition}
\begin{proof}
As with Proposition~\ref{psi_noncrossing_prop}, this is a straightforward consequence of Lemma~\ref{psi_right_nhbr_lem}. 
\end{proof}

We similarly get our desired conclusion. 

\begin{proposition}
\label{psi_ktermmin_prop}
The map $\psi$ restricts to a bijection between $k$-terminal-minimal diagrams of size $n$ and $(k-1)$-terminal-minimal diagrams of size $n-1$.
\end{proposition}

\begin{corollary}
The map $\psi^{k}$ on $k$-terminal diagrams restricts to a bijection between $k$-terminal-minimal diagrams of size $n+k$ and noncrossing diagrams of size $n$. In particular, $k$-terminal-minimal diagrams are counted by the Catalan numbers. 
\end{corollary}

While Lemma~\ref{topcycle_1term_lem} provides a structural characterization of 1-terminal-minimality, we have yet to obtain such a characterization for $k$-terminal-minimality. Our preliminary investigations indicate though that there should be a similar description in terms of $k$-terminal diagrams forbidding an infinite class of subdiagrams. 

\begin{conjecture}
There is a ``nice" infinite forbidden subdiagram characterization of $k$-terminal-minimality. 
\end{conjecture}

In 2015, Richard Stanley used ballot sequences and functions $f: \mathbb{N} \rightarrow \mathbb{N}$ satisfying $f(i) \leqslant i$ for all $i$ to define $2^{\aleph_{0}}$ chord diagram interpretations of the Catalan numbers $\{C_{n}\}_{n \geqslant 0}$. The noncrossing and nonnesting diagrams are recovered by setting $f(i) = 1$ and $f(i) = i$, respectively, so in a sense these interpretations lie between the two. Yet this construction only gives a finite number of interpretations for any fixed size $n$. With the set of $k$-terminal-minimal diagrams we have provided a countable infinity of combinatorial interpretations of $C_{n}$. Furthermore, combining Stanley's construction together with the map $\psi$ gives $2^{\aleph_{0}}$ combinatorial interpretations of the Catalan numbers with $\aleph_{0}$ interpretations at each fixed $n$. 

The above results characterize how $\psi$ acts on noncrossing diagrams; what about the other classical set of Catalan diagrams? The situation there is even nicer.

\begin{lemma} 
\label{psi_nesting_lem}
For $T \in \mathcal{T}_{n}$ and non-terminal $c, d \in T$, $d$ is nested under $c$ if and only if $\psi(d)$ is nested under $\psi(c)$.
\end{lemma} 
\begin{proof}
The equivalence follows immediately from the construction of $\psi$. 
\end{proof}

\begin{theorem} 
\label{psi_nest_cross_thm}
The map $\psi$ restricts to a bijection between 1-terminal diagrams with $n$ chords, $m$ crossings, and $k$ nestings and diagrams with $n-1$ chords, $m - n + 1$ crossings, and $k$ nestings.
\end{theorem}
\begin{proof}
This follows from Lemmas \ref{psi_right_nhbr_lem} and \ref{psi_nesting_lem} and the fact that the terminal chord of a 1-terminal diagram is not part of any nestings. 
\end{proof}

This result generalizes Theorem~\ref{psi_bijection_thm} and Propositions \ref{psi_noncrossing_prop} and \ref{psi_ktermmin_prop}. With it, one could begin to apply the large body of work (e.g. \cite{Chen2007, Pilaud2014}) on crossings and nestings in chord diagrams to 1-terminal diagrams. For example, Pilaud and Ru\'{e} \cite{Pilaud2014} used analytic combinatorial methods to derive an asymptotic estimate for the number of diagrams with $n$ chords and $m$ crossings; the above theorem can be used to extend this asymptotic estimate to 1-terminal diagrams with $n$ chords and $m$ crossings. 

\begin{corollary}
The map $\psi$ restricts to a bijection between connected nonnesting diagrams of size $n$ and nonnesting diagrams of size $n-1$.  
\end{corollary}
\begin{proof}
Applying Corollary~\ref{1term_char_cor} and Theorem~\ref{psi_nest_cross_thm} with $k = 0$ gives the result. 
\end{proof}

The above result uses the fact that connectivity is equivalent to 1-terminality for nonnesting diagrams. This actually points to a more general statement about nonnesting diagrams. 

\begin{lemma}
A nonnesting diagram is $k$-connected if and only if it is $k$-terminal and has at least $k$ chords. 
\end{lemma}
\begin{proof}
By Proposition~\ref{kterm_kconn_prop} it suffices to prove the ``only if" direction. Let $C$ be nonnesting and $k$-connected and $c$ be its root chord. By the Erd\"{o}s-Szekeres theorem the left and right neighborhoods of each chord of $C$ form cliques, implying by $k$-connectivity that the neighborhood of $c$ is a clique of size at least $k$. It follows that $c$ cannot be part of any minimal vertex cut, so $C - c$ is either complete or $k$-connected. In the former case we are done, while in the latter case we inductively get that $C - c$ is $k$-terminal, so $C$ is as well. 
\end{proof}

Combining this observation with Theorem~\ref{psi_nest_cross_thm} and Proposition~\ref{psi_kterm_prop}, we get the following. 

\begin{proposition}
The map $\psi$ restricts to a bijection between $k$-connected nonnesting diagrams of size $n$ and $(k-1)$-connected nonnesting diagrams of size $n-1$. 
\end{proposition}

\begin{corollary}
The map $\psi^{k}$ restricts to a bijection between $k$-connected nonnesting diagrams of size $n+k$ and nonnesting diagrams of size $n$. In particular, $k$-connected nonnesting diagrams are counted by the Catalan numbers. 
\end{corollary}






\subsection{Relationship with other double factorial objects}
\label{factorial_objects_sec}

In addition to chord diagrams, there are a number of classical combinatorial objects counted by the double factorials $(2n - 1)!!$. We will briefly focus on two of the most notable: increasing ordered trees and Stirling permutations. A (rooted) tree $T$ of size $n$ is {\em ordered} or {\em plane} if it is equipped with a total ordering of the children of each vertex $v \in T$. The tree $T$ is {\em increasing} if its vertices are labelled with the integers $0, 1, 2, \ldots, n-1$ such that the label of any child is greater than its parent. A {\em Stirling permutation} of size $n$ is a permutation of the multiset $\{1, 1, 2, 2, \ldots, n, n\}$ such that, for all $i$, all values between the two occurrences of $i$ are greater than $i$. Stirling permutations were introduced by Gessel and Stanley to give a combinatorial interpretation to the Stirling polynomials. They have since been studied extensively and generalized in multiple different ways (see e.g. \cite{Gessel1978, Bona2009, Janson2008, Janson2011}), as have increasing ordered trees (see e.g. \cite{Bergeron1992, Panholzer2007}). Write $\mathcal{S}_{n}$ and $\mathcal{I}_{n}$ for the sets of Stirling permutations and increasing ordered trees of size $n$, respectively. 

There are classic recursive constructions of both $\mathcal{C}_{n}$ and $\mathcal{S}_{n}$: insert a root chord and insert $nn$ into each element of $\mathcal{C}_{n-1}$ and $\mathcal{S}_{n-1}$, respectively, in all possible ways. Since there are always $2n - 1$ possible insertion places, this proves that both of these sets have cardinality $(2n - 1)!!$ and gives a recursively-defined bijection $\zeta$ between them. This leads into a simple characterization of 1-terminality and the bijection $\psi$ on Stirling permutations. 

\begin{proposition}
Diagram $C$ is 1-terminal if and only if $\zeta(C)$ has 1 as a prefix and suffix. In this case, $\zeta(\psi(C))$ is obtained from $\zeta(C)$ by removing both occurrences of 1 and normalizing the alphabet.  
\end{proposition}
\begin{proof}
Suppose $\zeta(C)$ does not have 1 as a prefix or suffix; by possibly taking the reverse, we may assume the former. Then, since $11$ clearly has 1 has a prefix, at some step in the iterative construction of $C$ a root chord is concatenated to the front of the currently placed chords, implying that it must be terminal in $C$ and not the terminal chord with greatest sink. So $C$ is not 1-terminal. On the other hand, if $C$ is 1-terminal, this never occurs, so it easily follows that 1 remains a prefix and suffix whenever a new element is placed and, therefore, is a prefix and suffix of $\zeta(C)$; this gives the first statement. For the second statement, let $C'$ be obtained by removing the root chord $c$ of $C$. Observe that $\psi(C') = \psi(C) - \psi(c)$ and $\zeta(C') = \zeta(C) - nn$. Then the result follows straightforwardly from induction. 
\end{proof}

The classic bijection $\eta$ between increasing ordered trees and Stirling permutations sends trees of size $n+1$ to permutations of size $n$: for a tree $T \in \mathcal{I}_{n+1}$, delete the root label and transfer the remaining labels from vertices to parent edges, then take a pre-order traversal around the tree, traversing each edge twice. The encountered sequence of labels is the Stirling permutation $\eta(T)$. 

As hinted at prior in the paper, and from the fact that there are $(2n-1)!!$ increasing ordered trees of size $n+1$, the decomposition from Section~\ref{bin_eq_sec} defined by the maps $\alpha$ and $\beta$ recursively defines a natural bijection $\theta$ between the set $\mathcal{T}_{n+1}$ of 1-terminal diagrams and the set $\mathcal{I}_{n+1}$ of increasing ordered trees. For each $C \in T_{n+1}$, apply $\alpha$ to get the 1-terminal subdiagrams $C_{1}, \ldots, C_{k}$ and intervals $I_{1}, \ldots, I_{k}$ of $[1, n]$. Then recursively apply $\theta$ to $C_{1}, \ldots, C_{k}$ to get increasing ordered trees $T_{1}, \ldots, T_{k}$, graft them to a new root of label 0 in that order, and finally apply an order-preserving bijection to reassign each subtree $T_{\ell}$ with the labels from $I_{\ell}$. The resulting increasing ordered tree $T$ is then set as the image of $C$ under $\theta$. We omit the proof that this actually defines a bijection between $\mathcal{T}_{n+1}$ and $\mathcal{I}_{n+1}$; it proceeds similarly to the proof of Theorem~\ref{diagram_triangle_thm}. Then, in this context, $\psi$ plays the role of extending this to a bijection to chord diagrams of size $n$. 

\begin{figure}[t]
\begin{center}
\includegraphics[scale=0.85]{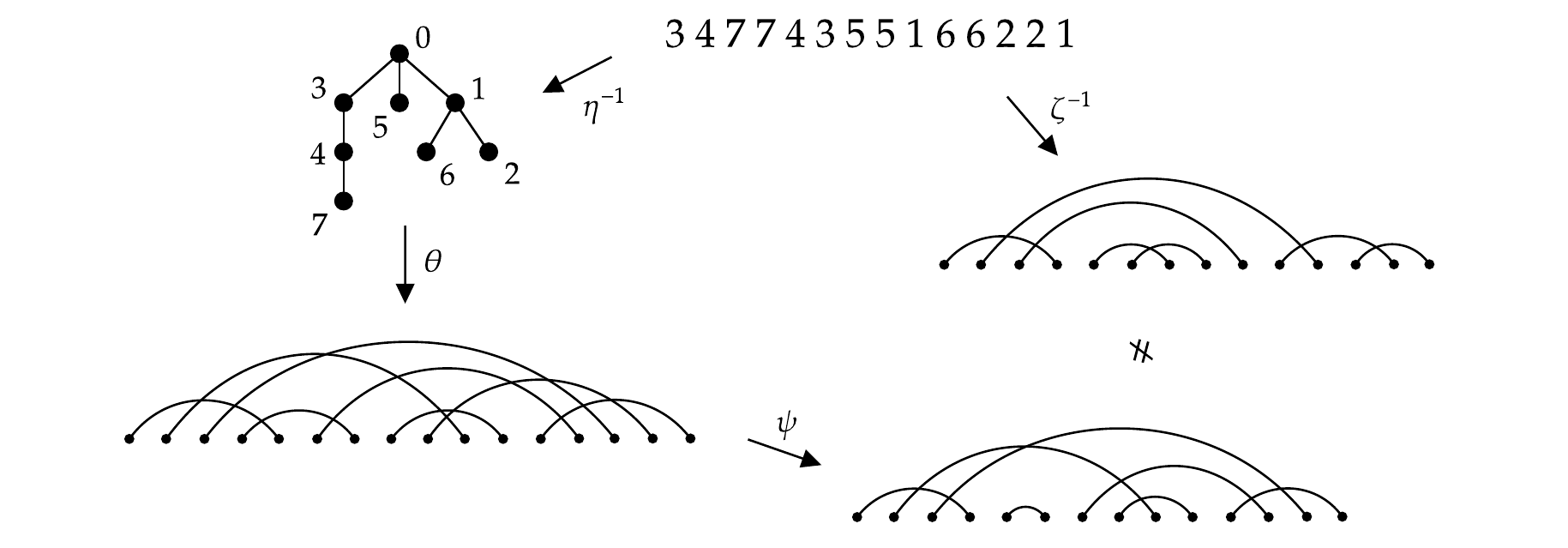}
\caption{A Stirling permutation $\sigma$ such that $\zeta^{-1}(\sigma) \neq (\psi \circ \theta \circ \eta^{-1})(\sigma)$.}
\label{double_factorial_maps_fig}
\end{center}
\end{figure}

From all of the above we get two bijections between Stirling permutations of size $n$ and chord diagrams with $n$ chords: $\zeta^{-1}$ and $\psi \circ \theta \circ \eta^{-1}$. These bijections are highly distinct, typically mapping a given Stirling permutation to two different diagrams; e.g. see Figure~\ref{double_factorial_maps_fig}. We can think of $\zeta^{-1}$ as encoding the {\em recursive view} of chord diagrams and $\psi \circ \theta \circ \eta^{-1}$ as encoding the {\em tree structure view} of chord diagrams. The map $\psi$ plays a role in both perspectives, forming a kind of bridge translating 1-terminality from the recursive view to the tree structure view. 

\subsection{Relationship with connectivity}
\label{psi_conn_sec}

We have seen how higher connectivity Catalan objects arise from $\psi$. Given the importance of connectedness to the story thus far, as well the appearance of other graph-theoretic properties and objects, it is natural to consider how connectivity might play a more general role. The notion is not new to chord diagrams. Daniel J. Kleitman \cite{Kleitman1970} argued that the number of $k$-connected diagrams is asymptotically $(2n-1)!!/e^{-k}$. For the $k = 2$ case Ali A. Mahmoud extended his result and used a decomposition of connected diagrams to obtain an asymptotic expansion for the number of 2-connected diagrams on $n$ chords. Here, we study the relationship between the map $\psi$ and connectivity. 

\begin{figure}[t]
\begin{center}
\includegraphics[scale=0.75]{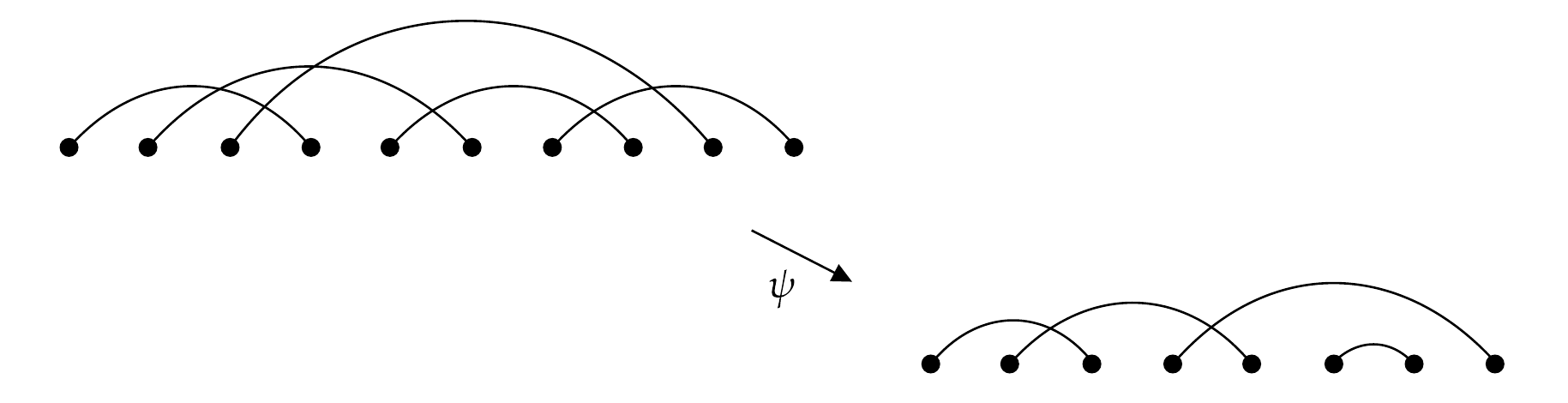}
\caption{An example of a 2-connected 1-terminal diagram mapped to a disconnected diagram by $\psi$.}
\label{psi_2conn_drop_fig}
\end{center}
\end{figure}

In light of what is known about $\psi$, especially Proposition \ref{kterm_kconn_prop} and \ref{psi_kterm_prop}, one might be initially tempted to conjecture that $\psi$ simply sends $k$-connected 1-terminal diagrams to $(k-1)$-connected diagrams, generalizing its behavior on $k$-terminal diagrams, but Figure~\ref{psi_2conn_drop_fig} gives a counterexample. So how much can the connectivity drop under $\psi$? This question, and the previous conjecture, was first put forth by Michael Borinsky.\footnote{Personal communication.} The following result answers this question succinctly. But, before stating the theorem, we need to define two types of diagrams for the proof. A {\em permutation diagram} is a chord diagram in which all the sinks proceed all the sources; in other words, a line can be drawn from the gap between two adjacent endpoints to immediately after the last sink that crosses every chord. A {\em shifted permutation diagram} is a 1-terminal diagram which becomes a permutation diagram after removing the terminal chord. By Lemma~\ref{order_agree_lem} a permutation diagram is 1-terminal if and only if the terminal chord crosses every other chord. Furthermore, it is easy to see that the neighborhoods of the last two chords of a shifted permutation diagram union to the entire diagram and that these two chords cross, and that these two properties determine whether a diagram is shifted permutation. 

\begin{figure}[t]
\begin{center}
\includegraphics[scale=0.75]{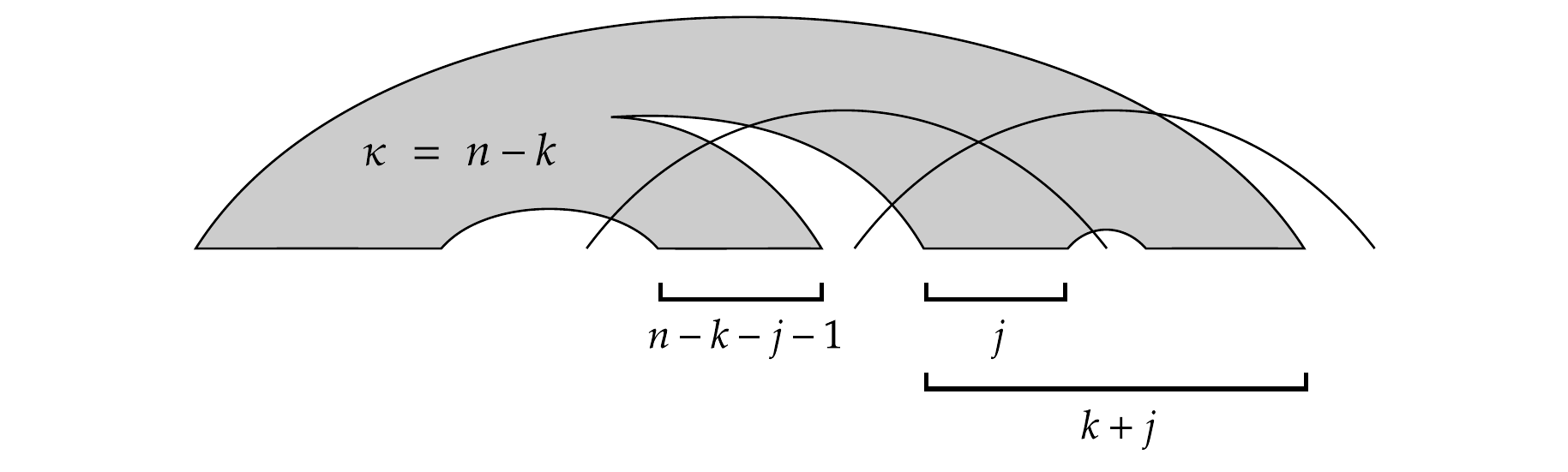}
\caption{A shifted permutation diagram of size $n$ and connectivity $n-k$ that maps under $\psi$ to a diagram of connectivity $j$. The first $n-2$ chords form an arbitrary $(n-k)$-connected permutation diagram.}
\label{psi_conn_example_fig}
\end{center} 
\end{figure}

\begin{theorem}
\label{psi_conn_thm}
The map $\psi$ sends each 1-terminal diagram of size $n$ and connectivity $n - k$ to a diagram with connectivity $n - 2k \leqslant j < n - k$, where $1 \leqslant k < n$. Furthermore, there is an example for each such $n$, $k$, and $j \geqslant 0$. 
\end{theorem}
\begin{proof}
Let $C$ be a 1-terminal diagram of size $n$ and connectivity $n-k$. By Lemma~\ref{psi_right_nhbr_lem}, the construction of $\psi$, and the fact that the terminal chord of $C$ has at least $n-k$ neighbors, the image $\psi(C)$ is obtained from $C$ by deleting one vertex and at most $k-1$ edges. Then clearly the connectivity of $\psi(C)$ is at least $\kappa(C) - k = n - 2k$, giving the lower bound. For the upper bound, we proceed via induction on $k$ and, for this purpose, also attach the statement that its inverse sends each 1-terminal diagram of size $n$ and connectivity $n-k$ to a diagram of size $n+1$ and connectivity at least $n-k+1$. Since $\psi$ is a bijection this is only a stronger statement in the context of the induction argument. Clearly the only $(n-1)$-connected diagram of size $n$, the complete diagram, maps to its size $n-1$ equivalent and vice versa under the inverse, as required. Now suppose to the contrary that $\kappa(\psi(C)) = n - k$. Since $n - k = n - 1 - (k - 1)$, the induction hypothesis implies that there exists an $(n-k+1)$-connected diagram $C'$ such that $\psi(C') = \psi(C)$. But $C'$ has strictly higher connectivity than $C$ so this contradicts the injectivity of $\psi$. 

For the second part of the result, consider a shifted permutation diagram $P$ with chords $c_{1} < c_{2} < \cdots < c_{n}$ such that $\{c_{1}, c_{2}, \ldots, c_{n-2}\}$ is a $(n-k)$-connected subdiagram, $|N(c_{n})| = k + j$, and $|N(c_{n-1}) \cap N(c_{n})| = j$. Clearly such a diagram exists; see Figure~\ref{psi_conn_example_fig}. Furthermore, $|N(c_{n-2})| = n - k$ and $k + j \geqslant n - k$, so $P$ has connectivity $n - k$ by Lemma~\ref{graph_conn_lem} with $A = \{c_{n-1}, c_{n}\}$. But $\psi(P)$ is obtained from $C$ by deleting $c_{n}$ and uncrossing $c_{n-1}$ with every chord in $N(c_{n-1}) - N(c_{n})$; in other words, removing the degree $j < n - k$ chord $c_{n-1}$ from $\psi(P)$ leaves an $(n - k)$-connected diagram. Thus $\psi(P)$ has connectivity $j$. 
\end{proof}

The nontrivial part of this theorem is that the connectivity necessarily decreases at all. If we view $\psi$ as an operation on graphs $G$ that deletes a vertex $v$ and one distinct edge incident to every vertex in $V(G) - N(v)$, then it is easy to see that such an operation does not always result in a drop in connectivity. The situation may change if we impose an order structure. As previously discussed, the directed intersection graph of a chord diagram is a directed acyclic graph. Every linear extension (e.g. the standard order, intersection order, peeling order) of the partial order determined by this directed acyclic graph is associated to a topological ordering of the graph. Then 1-terminal diagrams are the chord diagrams whose directed intersection graphs have exactly one sink, corresponding to the terminal chord. We can thus view $\psi$ as an operation on directed acyclic graphs with a single sink that, relative to a fixed topological ordering, deletes the last outgoing edge incident to each vertex and then removes the now isolated, previously unique sink vertex. For brevity, we will refer to such operations as {\em $\psi$-type maps} and such directed acyclic graphs as {\em 1-terminal}. 

\begin{question}
Does Theorem~\ref{psi_conn_thm} generalize to $\psi$-type maps on 1-terminal directed acyclic graphs?
\end{question}  

Most likely the answer to this question is no, but a counterexample eludes us. The proof of Theorem~\ref{psi_conn_thm} vitally uses the fact that $\psi$ is a bijection, but $\psi$-type maps are not bijections. In fact, $\psi$ itself is not even an bijection when viewed as acting on the directed intersection graph: $\{(1, 3), (2, 5), (4, 6)\}$ and $\{(1, 5), (2, 4), (3, 6)\}$ have non-isomorphic directed intersection graphs, but their images do not, and this example still works even when using the standard order to assign vertex labels and instead asking for distinct graphs. It is not clear how to prove Theorem~\ref{psi_conn_thm} without using the fact that $\psi$ is a bijection. 

Permutation diagrams are chord diagram representatives of permutation graphs, which have been extensively studied in the graph theory literature and play a central role in the study of circle graphs (e.g. \cite{Davies2020}). It is notable then that they also appear in this context and, in particular, are well-behaved with respect 1-terminality and $\psi$.

We end this section by asking whether it is possible to extend Theorem~\ref{psi_conn_thm} and obtain a useful characterization of the 1-terminal diagrams of size $n$ whose connectivity drops by a fixed integer $\ell$ under $\psi$. In particular, it would be of interest to enumerate these diagrams.




\section*{Acknowledgements}

The author would like to thank Karen Yeats and Nicholas Olson-Harris for valuable discussions on terminal chords, 1-terminal diagrams, the map $\psi$, and other topics in this paper. Parts of this work arose from these discussions. Karen Yeats also provided helpful comments on presentation.

\end{document}